\documentstyle[12pt]{article}
\input amssym.def
\begin{document}
\def \Z{\Bbb Z}
\def \C{\Bbb C}
\def \R{\Bbb R}
\def \Q{\Bbb Q}
\def \N{\Bbb N}
\def \wt{{\rm wt}}
\def \tr{{\rm tr}}
\def \span{{\rm span}}
\def \Res{{\rm Res}}
\def \End{{\rm End}\;}
\def \Ind {{\rm Ind}}
\def \Irr {{\rm Irr}}
\def \Aut{{\rm Aut}}
\def \Hom{{\rm Hom}}
\def \mod{{\rm mod}}
\def \ob{{\rm ob}\;}
\def \ann{{\rm Ann}}
\def \<{\langle} 
\def \>{\rangle} 
\def \t{\tau }
\def \a{\alpha }
\def \e{\epsilon }
\def \l{\lambda }
\def \L{\Lambda }
\def \g{\gamma}
\def \b{\beta }
\def \om{\omega }
\def \o{\omega }
\def \c{\chi}
\def \ch{\chi}
\def \cg{\chi_g}
\def \ag{\alpha_g}
\def \ah{\alpha_h}
\def \ph{\psi_h}

\def \bconj{\begin{conj}\label}
\def \econj{\end{conj}}
\def \be{\begin{equation}\label}
\def \ee{\end{equation}}
\def \bex{\begin{exa}\label}
\def \eex{\end{exa}}
\def \bl{\begin{lem}\label}
\def \el{\end{lem}}
\def \bt{\begin{thm}\label}
\def \et{\end{thm}}
\def \bp{\begin{prop}\label}
\def \ep{\end{prop}}
\def \br{\begin{rem}\label}
\def \er{\end{rem}}
\def \bc{\begin{coro}\label}
\def \ec{\end{coro}}
\def \bd{\begin{de}\label}
\def \ed{\end{de}}
\def \pf{{\bf Proof. }}
\def \voa{{vertex operator algebra}}

\newtheorem{thm}{Theorem}[section]
\newtheorem{prop}[thm]{Proposition}
\newtheorem{coro}[thm]{Corollary}
\newtheorem{conj}[thm]{Conjecture}
\newtheorem{exa}[thm]{Example}
\newtheorem{lem}[thm]{Lemma}
\newtheorem{rem}[thm]{Remark}
\newtheorem{de}[thm]{Definition}
\newtheorem{hy}[thm]{Hypothesis}
\makeatletter
\@addtoreset{equation}{section}
\def\theequation{\thesection.\arabic{equation}}
\makeatother
\makeatletter

\makeatletter
\newlength{\@pxlwd} \newlength{\@rulewd} \newlength{\@pxlht}
\catcode`.=\active \catcode`B=\active \catcode`:=\active \catcode`|=\active
\def\sprite#1(#2,#3)[#4,#5]{
   \edef\@sprbox{\expandafter\@cdr\string#1\@nil @box}
   \expandafter\newsavebox\csname\@sprbox\endcsname
   \edef#1{\expandafter\usebox\csname\@sprbox\endcsname}
   \expandafter\setbox\csname\@sprbox\endcsname =\hbox\bgroup
   \vbox\bgroup
      \catcode`.=\active\catcode`B=\active\catcode`:=\active\catcode`|=\active
      \@pxlwd=#4 \divide\@pxlwd by #3 \@rulewd=\@pxlwd
      \@pxlht=#5 \divide\@pxlht by #2
      \def .{\hskip \@pxlwd \ignorespaces}
\def B{\@ifnextchar B{\advance\@rulewd by \@pxlwd}{\vrule
         height \@pxlht width \@rulewd depth 0 pt \@rulewd=\@pxlwd}}
      \def :{\hbox\bgroup\vrule height \@pxlht width 0pt depth
0pt\ignorespaces}
      \def |{\vrule height \@pxlht width 0pt depth 0pt\egroup
         \prevdepth= -1000 pt}
   }
\def\endsprite{\egroup\egroup}
\catcode`.=12 \catcode`B=11 \catcode`:=12 \catcode`|=12\relax
\makeatother

\def\hboxtr{\FormOfHboxtr} 
\sprite{\FormOfHboxtr}(25,25)[0.5 em, 1.2 ex] 

:BBBBBBBBBBBBBBBBBBBBBBBBB |
:BB......................B |
:B.B.....................B |
:B..B....................B |
:B...B...................B |
:B....B..................B |
:B.....B.................B |
:B......B................B |
:B.......B...............B |
:B........B..............B |
:B.........B.............B |
:B..........B............B |
:B...........B...........B |
:B............B..........B |
:B.............B.........B |
:B..............B........B |
:B...............B.......B |
:B................B......B |
:B.................B.....B |
:B..................B....B |
:B...................B...B |
:B....................B..B |
:B.....................B.B |
:B......................BB |
:BBBBBBBBBBBBBBBBBBBBBBBBB |

\endsprite

\begin{center}{\Large \bf
Regular representations and Huang-Lepowsky tensor functors
for vertex operator algebras}
\end{center}
\begin{center}{Haisheng Li\footnote{Supported by NSF grant 
DMS-9970496 and a grant from Rutgers University Research Council}\\
Department of Mathematical Sciences\\
Rutgers University\\
Camden, NJ 08102}
\end{center}

\begin{abstract}
This is the second paper in a series to study regular representations
for vertex operator algebras.
In this paper, given a module $W$ for a vertex operator algebra $V$,
we construct, out of the dual space $W^{*}$,
a family of canonical (weak) $V\otimes V$-modules called ${\cal{D}}_{Q(z)}(W)$
parametrized by a nonzero complex number $z$.
We prove that for $V$-modules $W,W_{1}$
and $W_{2}$, a $Q(z)$-intertwining map of type
${W'\choose W_{1}W_{2}}$ in the sense of Huang and Lepowsky
exactly amounts to a 
$V\otimes V$-homomorphism from $W_{1}\otimes W_{2}$ to
${\cal{D}}_{Q(z)}(W)$ and that a $Q(z)$-tensor product
of $V$-modules $W_{1}$ and $W_{2}$ in the sense of Huang and Lepowsky
amounts to a universal from $W_{1}\otimes W_{2}$ to the functor
${\cal{F}}_{Q(z)}$, where ${\cal{F}}_{Q(z)}$ is a functor from 
the category of $V$-modules to
the category of weak $V\otimes V$-modules defined by
${\cal{F}}_{Q(z)}(W)={\cal{D}}_{Q(z)}(W')$ for a $V$-module $W$.
Furthermore, Huang-Lepowsky's $P(z)$ and $Q(z)$-tensor functors 
for the category of $V$-modules are extended 
to functors $T_{P(z)}$ and $T_{Q(z)}$ from 
the category of $V\otimes V$-modules to the category of $V$-modules.
It is proved that functors ${\cal{F}}_{P(z)}$ and ${\cal{F}}_{Q(z)}$ 
are right adjoints of $T_{P(z)}$ and $T_{Q(z)}$, 
respectively.
\end{abstract}

\section{Introduction}

Let $V$ be a vertex operator algebra. 
In [Li2], for a $V$-module $W$ and a nonzero complex number $z$
we constructed a weak $V\otimes V$-module called
${\cal{D}}_{P(z)}(W)$ out of the dual space $W^{*}$. We proved that
a $P(z)$-intertwining map of type ${W'\choose W_{1}W_{2}}$ 
in the sense of [H3, HL1-4], as proved in [HL2], which is 
the evaluation of an intertwining operator,
exactly amounts to a $V\otimes V$-homomorphism from $W_{1}\otimes W_{2}$
into ${\cal{D}}_{P(z)}(W)$ ($\subset W^{*}=\overline{W'}$).
We furthermore proved that for $V$-modules $W_{1}$ and $W_{2}$,
a $P(z)$-tensor product in the sense of [H3, HL1-4] exactly amounts to
a universal from $W_{1}\otimes W_{2}$ to the functor ${\cal{F}}_{P(z)}$,
where ${\cal{F}}_{P(z)}$ is a functor from the category of $V$-modules to
the category of $V\otimes V$-modules defined by
${\cal{F}}_{P(z)}(W)={\cal{D}}_{P(z)}(W')$ for a $V$-module $W$.

The objects $P(z)$ and $Q(z)$ are two special elements of the moduli space 
$K$ of spheres with (ordered) punctures  and 
local coordinates vanishing at those punctures (see [HL2]), 
where the three ordered punctures are $\infty, z,0$ for $P(z)$ 
and $z, \infty, ,0$ for $Q(z)$ with the standard local coordinates.
In the tensor product theory developed in [H3, HL1-4]
for a vertex operator algebra $V$,
$P(z)$ and $Q(z)$-tensor functors were constructed.
As remarked in [HL2], these different tensor functors will
play important roles in the formulation and 
constructions of the associativity and commutativity isomorphisms.

In this paper, a sequel to [Li2], we shall obtain parallel results
(with $P(z)$ being replaced by $Q(z)$);
we construct a weak $V\otimes V$-module called
${\cal{D}}_{Q(z)}(W)$ (out of $W^{*}$) such that
a $Q(z)$-intertwining map of type ${W'\choose W_{1}W_{2}}$
in the sense of Huang-Lepowsky exactly 
amounts to a $V\otimes V$-homomorphism from $W_{1}\otimes W_{2}$
into ${\cal{D}}_{Q(z)}(W)$ and such that a $Q(z)$-tensor product 
of $V$-modules $W_{1}$ and $W_{2}$ in the sense of Huang-Lepowsky 
exactly amounts to a universal from $W_{1}\otimes W_{2}$ to a 
functor ${\cal{F}}_{Q(z)}$. 

Just like in the tensor product theory ([H3],[HL1-5]), 
each of the $P(z)$ and $Q(z)$-notions 
has its own advantages in practice.
For example, the $Q(z)$-notions are
more natural in considering trace functions.
In fact, in the study of genus-zero correlation functions [Li3], 
it is more natural to use the notion of ${\cal{D}}_{Q(z)}(W)$.

As implied by the results of [Li2], the functor ${\cal{F}}_{P(z)}$
is essentially a right adjoint of Huang-Lepowsky $P(z)$-tensor 
functor. To make this statement precise, one needs to extend
Huang-Lepowsky $P(z)$-tensor functor to a functor from the category of
$V\otimes V$-modules to the category of $V$-modules.
(Note that the product category of the category of $V$-modules is a
natural subcategory of the category of $V\otimes V$-modules.)
In this paper, we carry out this exercise
to extend Huang-Lepowsky's $P(z)$ and $Q(z)$-tensor 
functors for the category of $V$-modules to functors $T_{P(z)}$ and $T_{Q(z)}$,
respectively, from the category ${\cal{C}}_{V\otimes V}^{o}$ to the category
of $V$-modules, where ${\cal{C}}_{V\otimes V}^{o}$ consists of
$V\otimes V$-modules on which $L(0)\otimes 1$ and $1\otimes L(0)$
act semisimply. It is proved that if $V$ is rational in the sense of [HL1-4],
functors ${\cal{F}}_{P(z)}$ and ${\cal{F}}_{Q(z)}$ are exactly right adjoints 
of $T_{P(z)}$ and $T_{Q(z)}$, respectively.
In establishing the functors $T_{P(z)}$ and $T_{Q(z)}$ (Theorems \ref{texthl}
and \ref{tQ(z)=P(z)}), instead of following [HL1-4] 
([HL2], Theorem 6.1 and [HL4], Theorem 13.9), we
use a different approach by employing the regular representations.
On the other hand,  presumably Theorems \ref{texthl}
and \ref{tQ(z)=P(z)} follow from [HL1-4] 
straightforward without any further efforts.
Furthermore, we have greatly taken the advantage of
[HL1-4] and we are motivated by their results 
for Propositions \ref{pp(z)q(z)} and \ref{pQ(z)=P(z)}
and Theorem \ref{tQ(z)=P(z)}.

As mentioned in [Li2], there is 
a classical analogue in the representation theory 
of Lie algebras for the left-right adjointness phenomenon.
(We thank James Lepowsky for bringing up this classical analogue issue.)
Let $\frak{g}$ be a Lie algebra and $U(\frak{g})$
the universal enveloping algebra. Usually,
$U(\frak{g})$ as a left (right) $U(\frak{g})$-module is called 
the left (right) regular representation of $\frak{g}$ (cf. [Di]).
Let $\Delta$ be the diagonal map from $\frak{g}$
to the product Lie algebra $\frak{g}\times \frak{g}$.
Since $\Delta$ is a Lie algebra homomorphism,
the diagonal map $\Delta$ gives rise 
to a natural functor $T$ from the category of 
$\frak{g}\times \frak{g}$-modules
to the category of $\frak{g}$-modules. 
Furthermore, the standard tensor functor of the category of $\frak{g}$-modules 
is the restriction of this functor $T$. On the other hand,
for any $\frak{g}$-module $U$, we have a coinduced module
$\Hom_{U(\frak{g})}(U(\frak{g}\times \frak{g}), U)$ (cf. [Di]).
Then the notion of coinduced module gives rise to a functor $E$ from 
the category of $\frak{g}$-modules to the category of 
$\frak{g}\times \frak{g}$-modules. From [Di] (Proposition 5.5.3),
the coinduced module functor $E$ is a right adjoint of the tensor functor $T$. 
(The usual induced module functor is a left adjoint of the tensor functor $T$.)

In the classical Lie theory, 
the Harish-Chandra bimodule theory is quite important.
In [FM], Frenkel and Malikov studied
Harish-Chandra categories by using Kazhdan-Lusztig tensor 
functor for affine Lie algebras [KL1-4] and obtained certain very
interesting results. (Kazhdan-Lusztig tensor functor 
and Huang-Lepowsky's tensor functors for the vertex operator algebras
associated to affine Lie algebras are closely related.)
In vertex operator algebra theory, analogues of many classical
fundamental notions and theories have been carried out.
On the other hand, due to the complexity of the new theory, certain
(straightforward) analogues do not make any sense.
A very interesting question is how much the classical 
Harish-Chandra bimodule theory can be carried out in the theory of
vertex operator algebras. We may address this issue in the future.

This paper is organized as follows: In Section 2, we recall
the construction of ${\cal{D}}_{P(z)}(W)$ and the main results from [Li2].
Sections 3 and 4 are devoted to the construction of ${\cal{D}}_{Q(z)}(W)$
and its relation with Huang-Lepowsky's notion of $Q(z)$-tensor product.
In Sections 5 and 6, Huang-Lepowsky's $Q(z)$ and $P(z)$-tensor functors
are extended to functors $T_{Q(z)}$ and $T_{P(z)}$ and the left-right
adjointness of functors are established.

\section{The weak $V\otimes V$-module ${\cal{D}}_{P(z)}(W)$}

In this section, we shall review the construction of the weak
$V\otimes V$-module ${\cal{D}}_{P(z)}(W)$ and the main results
obtained in [Li2]. This section does not contain any new result,
however we shall use the explicit construction in later sections.

As in [Li2], we use standard definitions and notations as given in [FLM, FHL].
A vertex operator algebra is denoted by $(V,Y,{\bf 1},\omega)$,
where ${\bf 1}$ is the vacuum vector and $\omega$ is the Virasoro element,
or simply by $V$.
We also use the notion of weak module as defined in [DLM]---A weak module
satisfies all the axioms defining the notion of a module 
given in [FLM] and [FHL] except those involving grading. Then
no grading is required for a weak module.

Since this paper deals with only the representation theory of vertex operator
algebras, {\em throughout this paper we fix a vertex operator algebra $V$.}

We typically use letters $x,y, x_{1},x_{2},\dots$ for mutually commuting
formal variables and $z,z_{0},\dots$ for complex numbers.
For a vector space $U$, $U[[x,x^{-1}]]$ is the vector space of all
(possibly doubly infinite) formal series with coefficients in $U$,
$U((x))$ is the space of formal Laurent series and
$U((x^{-1}))$ is the space of formal series truncated from above.
We shall also use the following standard 
formal variable convention (cf. [FLM]):
\begin{eqnarray}
& &(x_{1}-x_{2})^{n}=\sum_{i\ge 0}(-1)^{i}{n\choose i}x_{1}^{n-i}x_{2}^{i},\\
& &(x-z)^{n}=\sum_{i\ge 0}(-z)^{i}{n\choose i}x^{n-i},\\
& &(z-x)^{n}=\sum_{i\ge 0}(-1)^{i}z^{n-i}{n\choose i}x^{i}
\end{eqnarray}
for $n\in {\Z},\; z\in {\C}^{\times}$.

For vector spaces $U_{1},U_{2}$, a linear map $f\in \Hom(U_{1},U_{2})$ 
extends canonically to a linear map from $U_{1}[[x,x^{-1}]]$ to 
$U_{2}[[x,x^{-1}]]$. We shall use this canonical extension
without any comments.

Let $M$ be a vector space
and let $G$ be a linear map from $V$ to $(\End M)[[x,x^{-1}]]$.
For $v\in V$, we set (cf. [FHL, HL1])
\begin{eqnarray}\label{eoppositeoperator}
G^{o}(v,x)=G(e^{xL(1)}(-x^{-2})^{L(0)}v,x^{-1}),
\end{eqnarray}
which lies in $(\End M)[[x,x^{-1}]]$
because $e^{xL(1)}(-x^{-2})^{L(0)}v\in V[x,x^{-1}]$.
Then ([FHL], Proposition 5.3.1)
\begin{eqnarray}
(G^{o})^{o}=G.
\end{eqnarray}
In particular, if $(W,Y)$ is a weak $V$-module,  for $v\in V$,
\begin{eqnarray}
Y^{o}(v,x)=Y(e^{xL(1)}(-x^{-2})^{L(0)}v,x^{-1})\in (\End W)[[x,x^{-1}]].
\end{eqnarray}
As pointed out in [HL2], the proof Theorem 5.2.1 of [FHL] 
proves that $Y^{o}$ satisfies an opposite Jacobi identity (cf. (\ref{ehl})).
It is a viewpoint of [HL1] that 
the pair $(W,Y^{o})$ should be thought of as a right (weak) $V$-module. 

A {\em right } weak $V$-module is
a vector space $W$ equipped 
with a linear map $Y_{W}$ from $V$ to $(\End W)[[x,x^{-1}]]$
such that for $v\in V$, $w\in W$,
\begin{eqnarray}
& &Y_{W}(v,x)w\in W((x^{-1})),\\
& &Y_{W}({\bf 1},x)=1,
\end{eqnarray}
and such that the following {\em opposite Jacobi identity} holds for $u,v\in V$:
\begin{eqnarray}\label{ehl}
& &x_{0}^{-1}\delta\left(\frac{x_{1}-x_{2}}{x_{0}}\right)
Y_{W}(v,x_{2})Y_{W}(u,x_{1})
-x_{0}^{-1}\delta\left(\frac{x_{2}-x_{1}}{-x_{0}}\right)
Y_{W}(u,x_{1})Y_{W}(v,x_{2})\nonumber\\
&=&x_{2}^{-1}\delta\left(\frac{x_{1}-x_{0}}{x_{2}}\right)
Y_{W}(Y(u,x_{0})v,x_{2}).
\end{eqnarray}

The following result formulated in [Li2], due to 
[FHL, HL2], gives the equivalence between the notion of 
(left) $V$-module and the notion of right $V$-module:

\bp{pcontragredient}
Let $W$ be a vector space equipped with a linear map 
$Y_{W}$ from $V$ to $(\End W)[[x,x^{-1}]]$. 
Then
$(W,Y_{W})$ is a left (weak) $V$-module if and only if 
$(W,Y^{o}_{W})$ is a right (weak) $V$-module. 
\ep

Let $W$ be a weak $V$-module for now. Then
\begin{eqnarray}
Y^{o}(v,x)\in \Hom (W,W((x^{-1})))\;\;\;\mbox{ for }v\in V.
\end{eqnarray}
Furthermore, for any complex number $z_{0}$,
\begin{eqnarray}
Y^{o}(v,x+z_{0})\in\Hom (W,W((x^{-1}))),
\end{eqnarray}
where by definition
\begin{eqnarray}
Y^{o}(v,x+z_{0})w=(Y^{o}(v,y)w)|_{y=x+z_{0}}
\end{eqnarray}
for $w\in W$.
Let $U$ be a vector space, e.g., $U={\C}$.
For $v\in V,\;f\in \Hom(W,U),\; z_{0}\in \C$, the compositions
$fY^{o}(v,x)$ and $fY^{o}(v,x+z_{0})$
are elements of $(\Hom (W,U))[[x,x^{-1}]]$. 

We recall the following result from [Li2]:

\bp{ptranslationz}
Let $(W,Y_{W})$ be a right weak 
$V$-module and let $z_{0}$ be any complex number.
For $v\in V$, we set
\begin{eqnarray}
Y_{W}^{(z_{0})}(v,x)=Y_{W}(v,x+z_{0})=e^{z_{0}{d\over dx}}Y_{W}(v,x)\in 
(\End W)[[x,x^{-1}]].
\end{eqnarray}
Then the pair
$(W,Y_{W}^{(z_{0})})$ is also a right weak $V$-module.
\ep

Next we briefly recall the construction of the weak $V\otimes V$-module 
${\cal{D}}_{P(z)}(W)$ from [Li2].

\bd{drecall1}
{\em Let $W$ be a weak $V$-module
and $z$ a nonzero complex number. A {\em $P(z)$-linear functional}
on $W$ is a linear functional $\alpha$ on $W$ that satisfies the 
following condition: For each $v\in V$, there exist $k,l\in {\N}$ such that
\begin{eqnarray}
x^{l}(x-z)^{k}\<\alpha, Y^{o}(v,x)w\>\in {\C}[x]
\end{eqnarray}
for all $w\in W$, or what is equivalent to, the series
$\<\alpha, Y^{o}(v,x)w\>$, an element of ${\C}((x^{-1}))$,
 absolutely converges
in the domain $|x|>|z|$ to a rational function of the form 
\begin{eqnarray}
x^{-l}(x-z)^{-k}g(x),
\end{eqnarray}
where $g(x)\in {\C}[x]$.}
\ed

All $P(z)$-linear functionals on $W$ form a subspace of $W^{*}$, 
denoted by ${\cal{D}}_{P(z)}(W)$.

The following is an obvious characterization 
of a $P(z)$-linear functional
without involving matrix-coefficients.

\bl{lrecall1} {\rm [Li2]}
Let $W, z$ be given as before and let $\alpha\in W^{*}$. 
Then $\alpha\in {\cal{D}}_{P(z)}(W)$ if and only if
for $v\in V$, there exist $k,l\in {\N}$ such that
\begin{eqnarray}\label{ers1}
x^{l}(x-z)^{k}\alpha Y^{o}(v,x)\in W^{*}[[x]],
\end{eqnarray}
or equivalently, if and only if for $v\in V$, 
there exists $k\in {\N}$ such that
\begin{eqnarray}
(x-z)^{k}\alpha Y^{o}(v,x)\in W^{*}((x)).
\end{eqnarray}
\el

Denote by ${\C}(x)$ the algebra of rational functions of $x$.
The $\iota$-maps $\iota_{x;0}$ and $\iota_{x;\infty}$ from 
${\C}(x)$ to ${\C}[[x,x^{-1}]]$ were defined as follows:
For any rational function $f(x)$, 
$\iota_{x;0}f(x)$ is the Laurent series expansion of $f(x)$ at $x=0$
and $\iota_{x;\infty}f(x)$ is the Laurent series expansion 
of $f(x)$ at $x=\infty$. These are injective 
${\C}[x,x^{-1}]$-linear maps. Using
the formal variable convention, we have
\begin{eqnarray}
& &\iota_{x;0}\left((x-z)^{n}f(x)\right)=(-z+x)^{n}\iota_{x;0}f(x),\\
& &\iota_{x;\infty}\left((x-z)^{n}f(x)\right)=(x-z)^{n}\iota_{x;\infty}f(x)
\end{eqnarray}
for $n\in {\Z},\; z\in {\C}^{\times},\; f(x)\in {\C}(x)$.
Furthermore, for $z\in {\C},\; f(x)\in {\C}(x)$,
\begin{eqnarray}
& &\iota_{x;0}f(x)=(\iota_{y;\infty}f(y^{-1}))|_{y=x^{-1}},\label{ezeroinfty}\\
& &\iota_{x;\infty}f(x-z)=(\iota_{y;\infty}f(y))|_{y=x-z}.\label{einftyz}
\end{eqnarray}
Note that $\iota_{y;\infty}f(y)\in {\C}((y^{-1}))$, so that
$(\iota_{y;\infty}f(y))|_{y=x-z}$ exists in ${\C}[[x,x^{-1}]]$.

{}From the definition, 
for $\alpha\in {\cal{D}}_{P(z)}(W),\; v\in V,\; w\in W$,
$\<\alpha, Y^{o}(v,x)w\>$ lies in the range of $\iota_{x;\infty}$.
Then $\iota_{x;\infty}^{-1}\<\alpha, Y^{o}(v,x)w\>$ is a well 
defined element of ${\C}(x)$. It is also true that 
$\<\alpha, Y^{o}(v,x+z)w\>$ lies in the range of $\iota_{x;\infty}$. 

\bd{drecall2} {\rm [Li2]}
{\em For $v\in V$, $\alpha\in {\cal{D}}_{P(z)}(W)$, we define
$$Y_{P(z)}^{L}(v,x)\alpha,\;\;\;\; 
Y_{P(z)}^{R}(v,x)\alpha\in W^{*}[[x,x^{-1}]]$$
by
\begin{eqnarray}
\<Y_{P(z)}^{L}(v,x)\alpha,w\>&=&\iota_{x;0}\left(\iota_{x;\infty}^{-1}
\<\alpha,Y^{o}(v,x+z)w\>\right)\\
\<Y_{P(z)}^{R}(v,x)\alpha,w\>&=&\iota_{x;0}\iota_{x;\infty}^{-1}
\<\alpha,Y^{o}(v,x)w\>
\end{eqnarray}
for $w\in W$.}
\ed

We have ([Li2], Proposition 3.22):

\bp{prelationnew}
Let $v\in V,\; \alpha\in {\cal{D}}_{P(z)}(W)$. Then
\begin{eqnarray}\label{sjacnew}
& &x_{0}^{-1}\delta\left(\frac{x-z}{x_{0}}\right)
\alpha Y^{o}(v,x)-x_{0}^{-1}\delta\left(\frac{z-x}{-x_{0}}\right)
Y_{P(z)}^{R}(v,x)\alpha\nonumber\\
&=&z^{-1}\delta\left(\frac{x-x_{0}}{z}\right)Y_{P(z)}^{L}(v,x_{0})\alpha.
\end{eqnarray}
\ep

As an immediate consequence we have:

\bc{cbasic} {\rm [Li2]}
Let $v\in V,\; \alpha\in {\cal{D}}_{P(z)}(W)$. Then
\begin{eqnarray}
& &(-z+x)^{k}Y_{P(z)}^{R}(v,x)\alpha=(x-z)^{k}\alpha Y^{o}(v,x),\\
& &(z+x)^{l}Y_{P(z)}^{L}(v,x)\alpha=(x+z)^{l}\alpha Y^{o}(v,x+z),
\label{e2.26}
\end{eqnarray}
where $k$ and $l$ are any pair of integers such that (\ref{ers1}) holds.
(The integers $k$ and $l$ could be negative.)
\ec

We have ([Li2], Proposition 3.24):

\bp{precall2'} Let $W$ and $z$ be given as before. Then 
\begin{eqnarray}
Y_{P(z)}^{L}(v,x)\alpha,\;\;\;\; 
Y_{P(z)}^{R}(v,x)\alpha\in ({\cal{D}}_{P(z)}(W))((x))
\end{eqnarray}
for $v\in V,\; \alpha \in {\cal{D}}_{P(z)}(W)$. Furthermore,
\begin{eqnarray}
Y_{P(z)}^{L}(u,x_{1})Y_{P(z)}^{R}(v,x_{2})
=Y_{P(z)}^{R}(v,x_{2})Y_{P(z)}^{L}(u,x_{1})
\end{eqnarray}
on ${\cal{D}}_{P(z)}(W)$ for $u,v\in V$.
\ep

In view of Proposition \ref{precall2'}, $Y^{L}_{P(z)}$ and $Y^{R}_{P(z)}$
give rise to a well defined linear map 
\begin{eqnarray}
Y_{P(z)}=Y^{L}_{P(z)}\otimes Y^{R}_{P(z)}: V\otimes V\rightarrow 
\left(\End {\cal{D}}_{P(z)}(W)\right)[[x,x^{-1}]].
\end{eqnarray}
Then we have ([Li2], Theorem 3.17, 
Propositions 3.21 and 3.24 and Theorem 3.25):

\bt{trecall2} 
Let $W$ be a weak $V$-module and $z$ a nonzero complex number.
Then both the pairs $({\cal{D}}_{P(z)}(W), Y_{P(z)}^{L})$ and
$({\cal{D}}_{P(z)}(W), Y_{P(z)}^{R})$ carry 
the structure of a weak $V\otimes V$-module and
the pair $({\cal{D}}_{P(z)}(W), Y_{P(z)})$ 
carries the structure of a weak $V\otimes V$-module.
\et

As pointed out in [Li2], 
${\cal{D}}_{P(z)}: W\mapsto {\cal{D}}_{P(z)}(W)$
is a contravariant functor from the category of weak $V$-modules to
the category of weak $V\otimes V$-modules.
For a general $V$, even if $W$ is a $V$-module, ${\cal{D}}_{P(z)}(W)$
may be not an (ordinary) $V\otimes V$-module.
Recall from [Li2] that $R_{P(z)}(W)$ is the sum of 
$V\otimes V$-submodules of ${\cal{D}}_{P(z)}(W)$, on which
$L(0)\otimes 1$ and $1\otimes L(0)$ act semisimply.
Denote by ${\cal{C}}^{o}_{V\otimes V}$ the category of
$V\otimes V$-modules on which $L(0)\otimes 1$ and $1\otimes L(0)$ act
semisimply.

Vertex operator algebra $V$ is said to be {\em rational 
in the sense of Huang-Lepowsky} if $V$ has only finitely many
irreducible modules up to equivalence, every $V$-module is completely
reducible and if the fusion rule for every triple of irreducible
$V$-modules is finite.
If $V$ is rational in the sense of Huang-Lepowsky,
it was proved in [Li2] (Theorem 4.7) that for any $V$-module $W$, $R_{P(z)}(W)$ 
is a $V\otimes V$-module.
Then ${\cal{F}}_{P(z)}: W\mapsto R_{P(z)}(W')$ is a 
coinvariant functor from the category ${\cal{C}}_{V}$ of $V$-modules
to the category ${\cal{C}}^{o}_{V\otimes V}$.

\section{Weak $V\otimes V$-module ${\cal{D}}_{Q(z)}(W)$}

In this section, we shall construct a weak $V\otimes V$-module
${\cal{D}}_{Q(z)}(W)$ for any weak $V$-module $W$.

\bd{ddqzv} {\em  Let $W$ be a (left) weak $V$-module and $z$ a nonzero 
complex number. A linear functional $\alpha$ on $W$ is called 
{\em a $Q(z)$-linear functional} if for $v\in V,\; w\in W$,
the formal series
$$\<\alpha,Y(v,x)w\>\;\left(=\sum_{n\in \Z} \<\a,v_{n}w\> x^{-n-1}\right),$$
an element of ${\C}((x))$, absolutely converges in the domain
$0<|x|<|z|$ to a rational function 
$f_{\alpha}(v,w;x)\in {\C}[x,x^{-1},(x+z)^{-1}]$ such that
when $v$ is fixed with $w$ being free,
the orders of the possible poles at $x=-z, \infty$ 
of $f_{\alpha}(v,w;x)$ being viewed as meromorphic functions
on the sphere ${\C}\cup \{\infty\}$ are uniformally bounded.}
\ed

All $Q(z)$-linear functionals on $W$ clearly form a subspace of $W^{*}$,
which we denote by ${\cal{D}}_{Q(z)}(W)$.
Let $\alpha\in {\cal{D}}_{Q(z)}(W)$. From the definition, 
for $v\in V$, there exist
nonnegative integers $m$ and $l$ such that
\begin{eqnarray}\label{emly}
x^{-m}(x+z)^{l}\<\alpha, Y(v,x)w\>\in {\C}[x^{-1}]
\end{eqnarray}
for all $w\in W$. (Note that $m$ depends on $l$.)
Thus
\begin{eqnarray}\label{emlyo}
x^{-m}(x+z)^{l}\alpha Y(v,x)\in W^{*}[[x^{-1}]],
\end{eqnarray}
or equivalently,
\begin{eqnarray}\label{eqzformal}
(x+z)^{l}\alpha Y(v,x)\in W^{*}((x^{-1})).
\end{eqnarray}
Conversely, if $\alpha\in W^{*}$ satisfies the condition that
for every $v\in V$, there exists a nonnegative integer $l$
such that (\ref{eqzformal}) holds, then 
we easily see that $\alpha\in {\cal{D}}_{Q(z)}(W)$. 
(From (\ref{emlyo}) to (\ref{emly}) we also use the truncation 
condition: $Y(v,x)w\in W((x))$.)
This gives a slightly different characterization 
of a $Q(z)$-linear functional.

\bl{lqzequiv}
Let $\alpha\in W^{*}$. Then $\alpha\in {\cal{D}}_{Q(z)}(W)$ 
if and only if for every $v\in V$, there exists a nonnegative 
integer $l$ such that (\ref{eqzformal}) holds.$\;\;\;\;\Box$
\el

In the following we shall identify ${\cal{D}}_{Q(z)}(W)$
with ${\cal{D}}_{P(-z^{-1})}(W)$.

\bd{dytilde}
{\em Let $W$ and $z$ be given as before.
For $v\in V,\; \alpha\in {\cal{D}}_{Q(z)}(W)$, we define
\begin{eqnarray}
\tilde{Y}_{Q(z)}(v,x)\alpha\in W^{*}[[x,x^{-1}]]
\end{eqnarray}
by
\begin{eqnarray}\label{edefinition2.4}
\<\tilde{Y}_{Q(z)}(v,x)\alpha,w\>=\iota_{x;\infty}\iota_{x;0}^{-1}
\<\alpha,Y(v,x)w\>
\end{eqnarray}
for $w\in W$.}
\ed

Let $v\in V,\; \alpha\in {\cal{D}}_{Q(z)}(W)$ and let $l\in {\N}$ 
be such that (\ref{eqzformal}) holds.
Then
$$(x+z)^{l}\<\alpha,Y(v,x)w\>\in {\C}[x,x^{-1}]$$
for $w\in W$. Since $\iota$-maps are ${\C}[x,x^{-1}]$-linear and map $1$ to $1$, 
{}from (\ref{edefinition2.4}) we get
\begin{eqnarray}\label{etilderelation}
(x+z)^{l}\<\tilde{Y}_{Q(z)}(v,x)\alpha,w\>
=(x+z)^{l}\<\alpha,Y(v,x)w\>.
\end{eqnarray}
Thus
\begin{eqnarray}\label{etildeyay}
(x+z)^{l}\tilde{Y}_{Q(z)}(v,x)\alpha=(x+z)^{l}\alpha Y(v,x).
\end{eqnarray}

\bp{pp(z)q(z)}
Let $W$ and $z$ be given as before. Then
\begin{eqnarray}
{\cal{D}}_{Q(z)}(W)={\cal{D}}_{P(-z^{-1})}(W).
\end{eqnarray}
Furthermore,
\begin{eqnarray}\label{e3.9}
\tilde{Y}_{Q(z)}^{o}(v,x)\alpha=Y_{P(-z^{-1})}^{R}(v,x)\alpha
\end{eqnarray}
for $v\in V,\;\alpha\in {\cal{D}}_{Q(z)}(W)$.
\ep

\pf Let $\alpha\in W^{*}$. From Definition \ref{ddqzv},
$\alpha\in {\cal{D}}_{Q(z)}(W)$ 
if and only if for $v\in V,\;w\in W$, the formal series 
$$\<\alpha,Y(v,x^{-1})w\>$$
converges in the domain $|x|>|z^{-1}|$ to a rational function in 
${\C}[x,x^{-1},(x+z^{-1})^{-1}]$ such that
when $v$ is fixed with $w$ being free, the orders of possible poles at $0$ 
and $-z^{-1}$ are uniformally bounded.
Furthermore, since $e^{xL(1)}(-x^{-2})^{L(0)}v\in V[x,x^{-1}]$ for $v\in V$,
$\alpha\in {\cal{D}}_{Q(z)}(W)$ if and only if for $v\in V,\; w\in W$, 
the formal series $\<\alpha,Y^{o}(v,x)w\>$ absolutely converges 
in the domain $|x|>|z^{-1}|$ to a rational function in 
${\C}[x,x^{-1},(x+z^{-1})^{-1}]$ such that
when $v$ is fixed with $w$ being free, the orders of possible poles at $0$ 
and $-z^{-1}$ are uniformally bounded.
Thus, $\alpha$ is a $Q(z)$-linear functional
if and only if $\alpha$ is a $P(-z^{-1})$-linear functional. 
This proves the first part of the proposition. 

For the second part, from (\ref{edefinition2.4}) 
(also using (\ref{ezeroinfty})) we have
\begin{eqnarray}
\<\tilde{Y}_{Q(z)}(v,x^{-1})\alpha,w\>
=\iota_{x;0}\iota_{x;\infty}^{-1}\<\alpha,Y(v,x^{-1})w\>.
\end{eqnarray}
Since $e^{xL(1)}(-x^{-2})^{L(0)}v\in V[x,x^{-1}]$ for $v\in V$, 
and the $\iota$-maps are ${\C}[x,x^{-1}]$-linear and send $1$ to $1$,
we have
\begin{eqnarray}
\<\tilde{Y}_{Q(z)}^{o}(v,x)\alpha,w\>
=\iota_{x;0}\iota_{x;\infty}^{-1}\<\alpha,Y^{o}(v,x)w\>.
\end{eqnarray}
Then (\ref{e3.9}) follows immediately from the definition 
of $Y_{P(-z^{-1})}^{R}$.
$\;\;\;\;\Box$

Because 
$({\cal{D}}_{P(-z^{-1})}(W),Y_{P(-z^{-1})}^{R})$
carries the structure of a left weak 
$V$-module (Theorem \ref{trecall2}), 
in view of Propositions \ref{pp(z)q(z)} and \ref{pcontragredient}
we immediately have (cf. [Li2], Theorem 3.17):

\bt{tqzleft}
Let $W$ be a (left) weak $V$-module and let $z$ be a nonzero 
complex number. Then 
the pair $({\cal{D}}_{Q(z)}(W),\tilde{Y}_{Q(z)})$
carries the structure of a right weak $V$-module. $\;\;\;\;\Box$
\et

Next, we are going to use $\tilde{Y}_{Q(z)}$ to define two left 
module structures on ${\cal{D}}_{Q(z)}(W)$.

\bd{dyqzlr}
{\em For $v\in V,\;\alpha\in {\cal{D}}_{Q(z)}(W)$,
we define 
\begin{eqnarray}
Y_{Q(z)}^{R}(v,x)\alpha,\;\;Y_{Q(z)}^{L}(v,x)\alpha\in W^{*}[[x,x^{-1}]]
\end{eqnarray}
by
\begin{eqnarray}
& &\<(Y^{L}_{Q(z)})^{o}(v,x)\alpha,w\>=\iota_{x;\infty}
\left(\iota_{y;0}^{-1}\<\alpha,Y(v,y)w\>\right)|_{y=x-z},\label{edefqzl}\\
& &\<Y^{R}_{Q(z)}(v,x)\alpha,w\>=\iota_{x;0}
\left(\iota_{y;0}^{-1}\<\alpha,Y(v,y)w\>\right)|_{y=x-z}\label{edefqzr}
\end{eqnarray}
for $w\in W$. }
\ed

Let $v\in V,\; \alpha\in {\cal{D}}_{Q(z)}(W),\;w\in W$. 
Using (\ref{edefqzl}) and (\ref{edefinition2.4}), we get
\begin{eqnarray}
\<(Y^{L}_{Q(z)})^{o}(v,x)\alpha,w\>
&=&\left(\iota_{y;\infty}
\iota_{y;0}^{-1}\<\alpha,Y(v,y)w\>\right)|_{y=x-z}\nonumber\\
&=&\left(\<\tilde{Y}_{Q(z)}(v,y)\alpha,w\>\right)|_{y=x-z}\nonumber\\
&=&\<\tilde{Y}_{Q(z)}(v,x-z)\alpha,w\>\nonumber\\
&=&\<\tilde{Y}_{Q(z)}^{(-z)}(v,x)\alpha,w\>.
\end{eqnarray}
Notice that $\<\tilde{Y}_{Q(z)}(v,x-z)\alpha,w\>$ exists
because $\tilde{Y}_{Q(z)}(v,y)\alpha\in {\cal{D}}_{Q(z)}(W)((y^{-1}))$.
(But $\<\alpha,Y(v,x-z)w\>$ does not exist.)
Thus
\begin{eqnarray}\label{etildeqz=left}
(Y^{L}_{Q(z)})^{o}(v,x)\alpha=\tilde{Y}_{Q(z)}^{(-z)}(v,x)\alpha.
\end{eqnarray}
Similarly, using (\ref{edefqzr}) and (\ref{edefinition2.4})
(also using (\ref{einftyz})) we get
\begin{eqnarray}
\<Y^{R}_{Q(z)}(v,x)\alpha,w\>
&=&\iota_{x;0}\left(\iota_{y;\infty}^{-1}
\<\tilde{Y}_{Q(z)}(v,y)\alpha,w\>\right)|_{y=x-z}\nonumber\\
&=&\iota_{x;0}\iota_{x;\infty}^{-1}
\<\tilde{Y}_{Q(z)}(v,x-z)\alpha,w\>\nonumber\\
&=&\iota_{x;0}\iota_{x;\infty}^{-1}
\<\tilde{Y}_{Q(z)}^{(-z)}(v,x)\alpha,w\>.
\end{eqnarray}
That is,
\begin{eqnarray}\label{eqzrighttilde}
\iota_{x;0}^{-1}\<Y^{R}_{Q(z)}(v,x)\alpha,w\>
=\iota_{x;\infty}^{-1}\<\tilde{Y}_{Q(z)}^{(-z)}(v,x)\alpha,w\>.
\end{eqnarray}
Notice that
\begin{eqnarray}\label{eqztilde=Yv}
\iota_{x;0}^{-1}\<\alpha,Y(v,x)w\>
=\iota_{x;\infty}^{-1}\<\tilde{Y}_{Q(z)}(v,x)\alpha,w\>
=\left(\iota_{y;\infty}^{-1}
\<\tilde{Y}_{Q(z)}^{(-z)}(v,y)\alpha,w\>\right)|_{y=x+z}.
\end{eqnarray}
Combining (\ref{eqzrighttilde}) with (\ref{eqztilde=Yv}), then
using the fundamental properties of the delta function 
([FHL], Proposition 3.1.1.) we get
\begin{eqnarray}\label{eqzrelation-1withw}
& &x_{0}^{-1}\delta\left(\frac{x_{1}-z}{x_{0}}\right)
\<\tilde{Y}_{Q(z)}^{(-z)}(v,x_{1})\alpha,w\>
-x_{0}^{-1}\delta\left(\frac{z-x_{1}}{-x_{0}}\right)
\<Y_{Q(z)}^{R}(v,x_{1})\alpha,w\>\nonumber\\
&=&z^{-1}\delta\left(\frac{x_{1}-x_{0}}{z}\right)
\<\alpha, Y(v,x_{0})w\>.
\end{eqnarray}
By using the delta-function substitution property,
dropping out $w$, we get another version:
\begin{eqnarray}\label{eqzrelation-1}
& &x_{0}^{-1}\delta\left(\frac{x_{1}-z}{x_{0}}\right)
\tilde{Y}_{Q(z)}(v,x_{0})\alpha
-x_{0}^{-1}\delta\left(\frac{z-x_{1}}{-x_{0}}\right)
Y_{Q(z)}^{R}(v,x_{1})\alpha\nonumber\\
&=&z^{-1}\delta\left(\frac{x_{1}-x_{0}}{z}\right)
\alpha Y(v,x_{0}).
\end{eqnarray}
Then we immediately have:

\bp{pypzlr}
For $v\in V,\; \alpha\in {\cal{D}}_{Q(z)}(W)$,
\begin{eqnarray}\label{eqzrelation}
& &z^{-1}\delta\left(\frac{x_{1}-x_{0}}{z}\right)\alpha Y(v,x_{0})
\nonumber\\
&=&x_{0}^{-1}\delta\left(\frac{x_{1}-z}{x_{0}}\right)
(Y_{Q(z)}^{L})^{o}(v,x_{1})\alpha
-x_{0}^{-1}\delta\left(\frac{z-x_{1}}{-x_{0}}\right)
Y_{Q(z)}^{R}(v,x_{1})\alpha. \;\;\;\;\Box
\end{eqnarray}
\ep

\bp{pqzleft} Let $W$ be a (left) weak $V$-module and 
let $z$ be a nonzero complex number. 
Then the pair $({\cal{D}}_{Q(z)}(W), Y_{Q(z)}^{L})$ 
carries the structure of a (left) weak $V$-module.
\ep

\pf Because $\tilde{Y}_{Q(z)}$ gives rise to a right weak 
$V$-module structure on ${\cal{D}}_{Q(z)}(W)$ (Theorem \ref{tqzleft}), 
by Proposition \ref{ptranslationz},
$\tilde{Y}_{Q(z)}^{(-z)}$ also gives rise to a right 
weak $V$-module structure on ${\cal{D}}_{Q(z)}(W)$.
Furthermore, by Proposition \ref{pcontragredient}, $Y_{Q(z)}^{L}$ gives 
rise to a left 
weak $V$-module structure on ${\cal{D}}_{Q(z)}(W)$ because
$(Y_{Q(z)}^{L})^{o}=\tilde{Y}_{Q(z)}^{(-z)}$. $\;\;\;\;\Box$

We shall prove that the pair $({\cal{D}}_{Q(z)}(W),Y^{R}_{Q(z)})$ 
also carries the structure of a left weak $V$-module. First we have:

\bl{lqzrighttruncation}
Let $W,z$ be given as before. Then 
for $v\in V,\;\alpha\in {\cal{D}}_{Q(z)}(W),\;n\in {\Z}$,
\begin{eqnarray}
& &(z+x)^{n}\alpha Y(v,x)\in {\cal{D}}_{Q(z)}(W)[[x,x^{-1}]],\label{ealphayv}\\
& &Y^{R}_{Q(z)}(v,x)\alpha\in {\cal{D}}_{Q(z)}(W)((x)).\label{eqzyrtruncation}
\end{eqnarray}
\el

\pf Let $u\in V,\; m\in {\Z}$. From Borcherds' commutator formula we have
\begin{eqnarray}\label{ecommutator}
\alpha v_{m}Y(u,x)=\alpha Y(u,x)v_{m}-
\sum_{i\ge 0}{m\choose i}x^{m-i}\alpha Y(v_{i}u,x).
\end{eqnarray}
Since $v_{i}u\ne 0$ for only finitely many $i\ge 0$, 
there exists $l\in {\N}$ such that
\begin{eqnarray}
(x+z)^{l}\alpha Y(u,x),\;\; (x+z)^{l}\alpha Y(v_{i}u,x)\in W^{*}((x^{-1}))
\end{eqnarray}
for all $i\ge 0$. Then using (\ref{ecommutator}) we get
\begin{eqnarray}\label{ealphavmin}
(x+z)^{l}\alpha v_{m}Y(u,x)\in W^{*}((x^{-1}))
\end{eqnarray}
for all $m\in {\Z}$. Thus $\alpha v_{m}\in {\cal{D}}_{Q(z)}(W)$.
Notice that for $n\in {\Z}$, 
$$(z+x_{0})^{n}\alpha Y(v,x_{0})\;\;\mbox{ exists}$$
in $W^{*}[[x_{0},x_{0}^{-1}]]$ though its coefficient of $x_{0}^{t}$
for $t\in {\Z}$ is generally an infinite sum of $\alpha v_{m}$'s.
Furthermore, from (\ref{ealphavmin}) we have
\begin{eqnarray}
(x+z)^{l}(z+x_{0})^{n}\alpha Y(v,x_{0}) Y(u,x)
\in (W^{*}[[x_{0},x_{0}^{-1}]])((x^{-1})).
\end{eqnarray}
Then (\ref{ealphayv}) follows from Lemma \ref{lqzequiv}.

Let $l\in {\N}$ be such that (\ref{etildeyay}) holds.
Taking $\Res_{x_{0}}$ from (\ref{eqzrelation-1}),
using the fundamental properties of delta functions,
and using (\ref{etildeyay}) and (\ref{ealphayv}) we get
\begin{eqnarray}\label{eqzyrtruncation2}
& &Y^{R}_{Q(z)}(v,x_{1})\alpha\nonumber\\
&=&\Res_{x_{0}}x_{1}^{-1}\delta\left(\frac{x_{0}+z}{x_{1}}\right)
\tilde{Y}_{Q(z)}(v,x_{0})\alpha
-\Res_{x_{0}}x_{1}^{-1}\delta\left(\frac{z+x_{0}}{x_{1}}\right)
\alpha Y(v,x_{0})\nonumber\\
&=&\Res_{x_{0}}\sum_{n\in {\Z}}x_{1}^{-n-1}\left(
(x_{0}+z)^{n}\tilde{Y}_{Q(z)}(v,x_{0})\alpha
-(z+x_{0})^{n}\alpha Y(v,x_{0})\right)\nonumber\\
&=&\Res_{x_{0}}\sum_{n<l}x_{1}^{-n-1}\left(
(x_{0}+z)^{n}\tilde{Y}_{Q(z)}(v,x_{0})\alpha
-(z+x_{0})^{n}\alpha Y(v,x_{0})\right).
\end{eqnarray}
Noticing that for any $n\in {\Z}$, every coefficient of monomials in $x$
of $(x+z)^{n}\tilde{Y}_{Q(z)}(v,x)\alpha$ is a finite sum of the coefficients
of $\tilde{Y}_{Q(z)}(v,x)\alpha$, with Theorem \ref{tqzleft} we get
$$(x+z)^{n}\tilde{Y}_{Q(z)}(v,x)\alpha
\in {\cal{D}}_{Q(z)}(W)((x^{-1})).$$
Then (\ref{eqzyrtruncation}) immediately follows from (\ref{ealphayv}) and  
(\ref{eqzyrtruncation2}).
$\;\;\;\;\Box$

We have the following commutativity relations:

\bp{pqzleftrightcomm}
Let $W$  be a (left) weak $V$-module and let $z$ be a nonzero complex number. 
Then
\begin{eqnarray}
& &\tilde{Y}_{Q(z)}(u,x_{1})Y_{Q(z)}^{R}(v,x_{2})
=Y_{Q(z)}^{R}(v,x_{2})\tilde{Y}_{Q(z)}(u,x_{1}),\label{etildeyyrcomm}\\
& &Y_{Q(z)}^{L}(u,x_{1})Y_{Q(z)}^{R}(v,x_{2})
=Y_{Q(z)}^{R}(v,x_{2})Y_{Q(z)}^{L}(u,x_{1})\label{eqzylyrcomm}
\end{eqnarray}
on ${\cal{D}}_{Q(z)}(W)$ for $u,v\in V$.
\ep

\pf In view of the connection between $Y_{Q(z)}^{L}(u,x)$ and 
$\tilde{Y}_{Q(z)}(v,x)$ (recall (\ref{etildeqz=left})), 
we only need to prove (\ref{etildeyyrcomm}).
The proof is similar to that of Proposition 3.24 of [Li2].

It follows (cf. [DL]) from the (opposite) Jacobi identity that
there exists $k\in {\N}$ such that
\begin{eqnarray}
& &(x_{0}-x_{2})^{k}\alpha Y(u,x_{2})Y(v,x_{0})
=(x_{0}-x_{2})^{k}\alpha Y(v,x_{0})Y(u,x_{2}),\\
& &(x_{0}-x_{2})^{k}
\tilde{Y}_{Q(z)}(u,x_{2})\tilde{Y}_{Q(z)}(v,x_{0})\alpha
=(x_{0}-x_{2})^{k}\tilde{Y}_{Q(z)}(v,x_{0})\tilde{Y}_{Q(z)}(u,x_{2})\alpha.
\hspace{1cm}
\end{eqnarray}
In view of (\ref{ealphavmin}) and (\ref{etildeyay}) there exists $l\in {\N}$
such that
\begin{eqnarray}
& &(x_{2}+z)^{l}\tilde{Y}_{Q(z)}(u,x_{2})\alpha
=(x_{2}+z)^{l}\alpha Y(u,x_{2}),\\
& &(x_{2}+z)^{l}\tilde{Y}_{Q(z)}(u,x_{2})(\alpha Y(v,x_{0}))
=(x_{2}+z)^{l}\alpha Y(v,x_{0})Y(u,x_{2}).
\end{eqnarray}
Using (\ref{eqzrelation-1}) together with all the identities above we have
\begin{eqnarray}
& &x_{0}^{-1}\delta\left(\frac{z-x_{1}}{-x_{0}}\right)
(x_{0}-x_{2})^{k}(x_{2}+z)^{l}
\tilde{Y}_{Q(z)}(u,x_{2})Y_{Q(z)}^{R}(v,x_{1})\alpha\nonumber\\
&=&x_{0}^{-1}\delta\left(\frac{x_{1}-z}{x_{0}}\right)
(x_{0}-x_{2})^{k}(x_{2}+z)^{l}
\tilde{Y}_{Q(z)}(u,x_{2})\tilde{Y}_{Q(z)}(v,x_{0})\alpha\nonumber\\
& &-z^{-1}\delta\left(\frac{x_{1}-x_{0}}{z}\right)
(x_{0}-x_{2})^{k}(x_{2}+z)^{l}
\tilde{Y}_{Q(z)}(u,x_{2})(\alpha Y(v,x_{0}))\nonumber\\
&=&x_{0}^{-1}\delta\left(\frac{x_{1}-z}{x_{0}}\right)
(x_{0}-x_{2})^{k}(x_{2}+z)^{l}
\tilde{Y}_{Q(z)}(u,x_{2})\tilde{Y}_{Q(z)}(v,x_{0})\alpha\nonumber\\
& &-z^{-1}\delta\left(\frac{x_{1}-x_{0}}{z}\right)
(x_{0}-x_{2})^{k}(x_{2}+z)^{l}
\alpha Y(v,x_{0})Y(u,x_{2})\nonumber\\
&=&x_{0}^{-1}\delta\left(\frac{x_{1}-z}{x_{0}}\right)
(x_{0}-x_{2})^{k}(x_{2}+z)^{l}
\tilde{Y}_{Q(z)}(v,x_{0})\tilde{Y}_{Q(z)}(u,x_{2})\alpha\nonumber\\
& &-z^{-1}\delta\left(\frac{x_{1}-x_{0}}{z}\right)
(x_{0}-x_{2})^{k}(x_{2}+z)^{l}
\alpha Y(u,x_{2})Y(v,x_{0})\nonumber\\
&=&x_{0}^{-1}\delta\left(\frac{x_{1}-z}{x_{0}}\right)
(x_{0}-x_{2})^{k}(x_{2}+z)^{l}
\tilde{Y}_{Q(z)}(v,x_{0})\tilde{Y}_{Q(z)}(u,x_{2})\alpha\nonumber\\
& &-z^{-1}\delta\left(\frac{x_{1}-x_{0}}{z}\right)
(x_{0}-x_{2})^{k}(x_{2}+z)^{l}
(\tilde{Y}_{Q(z)}(u,x_{2})\alpha) Y(v,x_{0})\nonumber\\
&=&x_{0}^{-1}\delta\left(\frac{z-x_{1}}{-x_{0}}\right)
(x_{0}-x_{2})^{k}(x_{2}+z)^{l}
Y_{Q(z)}^{R}(v,x_{1})\tilde{Y}_{Q(z)}(u,x_{2})\alpha.
\end{eqnarray}
By taking $\Res_{x_{0}}$ and using the fundamental properties 
of delta functions we get
\begin{eqnarray}
& &(-z+x_{1}-x_{2})^{k}(x_{2}+z)^{l}
\tilde{Y}_{Q(z)}(u,x_{2})Y_{Q(z)}^{R}(v,x_{1})\alpha\nonumber\\
&=&(-z+x_{1}-x_{2})^{k}(x_{2}+z)^{l}
Y_{Q(z)}^{R}(v,x_{1})\tilde{Y}_{Q(z)}(u,x_{2})\alpha.
\end{eqnarray}
Then multiplying both sides by $(-z+x_{1}-x_{2})^{-k}(x_{2}+z)^{-l}$
we obtain (\ref{etildeyyrcomm}), noting that we are allowed to  multiply. $\;\;\;\;\Box$

Now we are ready to prove:

\bp{pqzright} Let $W$ be a (left) weak $V$-module and 
let $z$ be a nonzero complex number. Then
the pair $({\cal{D}}_{Q(z)}(W), Y_{Q(z)}^{R})$ 
carries the structure of a (left) weak $V$-module.
\ep

\pf It is clear that $Y_{Q(z)}^{R}({\bf 1},x)=1$ because 
$\tilde{Y}_{Q(z)}({\bf 1},x)=1$. With Lemma \ref{lqzrighttruncation},
we only need to prove the Jacobi identity. 

For convenience, let us locally use $\overline{Y}$ 
for $\tilde{Y}_{Q(z)}^{(-z)}$. That is, 
$({\cal{D}}_{Q(z)}(W),\overline{Y})$ is a right weak $V$-module. 
Furthermore, for $v\in V,\; w\in W$, if $l$ is a nonnegative integer 
such that $x^{l}Y(v,x)w\in W[[x]]$, then by taking 
$\Res_{x_{0}}x_{0}^{l}$ from (\ref{eqzrelation-1withw}) we get
\begin{eqnarray}\label{eyqzr=y}
(x-z)^{l}\<Y_{Q(z)}^{R}(v,x)\alpha,w\>=(x-z)^{l}\<\overline{Y}(v,x)\alpha,w\>
\end{eqnarray}
for {\em all} $\alpha\in {\cal{D}}_{Q(z)}(W)$. Furthermore, from (\ref{etildeyyrcomm})
$\bar{Y}$ and $R_{Q(z)}^{R}$ commute.

Now, let $u,v\in V,\; \alpha\in {\cal{D}}_{Q(z)}(W),\; w\in W$.
Let $l\in {\N}$ be such that
\begin{eqnarray}
x^{l}Y(u,x)w,\;\; x^{l}Y(v,x)w\in W[[x]].
\end{eqnarray}
Then
\begin{eqnarray}
& &(x_{1}-z)^{l}
\<Y^{R}_{Q(z)}(u,x_{1})Y^{R}_{Q(z)}(v,x_{2})\alpha,w\>
=(x_{1}-z)^{l}
\<\bar{Y}(u,x_{1})Y^{R}_{Q(z)}(v,x_{2})\alpha,w\>,\hspace{1cm}\\
& &(x_{2}-z)^{l}\<Y^{R}_{Q(z)}(v,x_{2})\tilde{Y}(u,x_{1})\alpha,w\>
=(x_{2}-z)^{l}\<\bar{Y}(v,x_{2})\tilde{Y}(u,x_{1})\alpha,w\>
\end{eqnarray}
because 
$$Y^{R}_{Q(z)}(v,x)\alpha,\;\;\bar{Y}(u,x)\alpha\in 
{\cal{D}}_{Q(z)}(W)[[x,x^{-1}]].$$
Then 
\begin{eqnarray}\label{eqzrightfirst}
& &(x_{1}-z)^{l}(x_{2}-z)^{l}
\<Y^{R}_{Q(z)}(u,x_{1})Y^{R}_{Q(z)}(v,x_{2})\alpha,w\>\nonumber\\
&=&(x_{1}-z)^{l}(x_{2}-z)^{l}
\<\bar{Y}(u,x_{1})Y^{R}_{Q(z)}(v,x_{2})\alpha,w\>\nonumber\\
&=&(x_{1}-z)^{l}(x_{2}-z)^{l}
\<Y^{R}_{Q(z)}(v,x_{2})\bar{Y}(u,x_{1})\alpha,w\>\nonumber\\
&=&(x_{1}-z)^{l}(x_{2}-z)^{l}
\<\bar{Y}(v,x_{2})\bar{Y}(u,x_{1})\alpha,w\>.
\end{eqnarray}
Similarly we have
\begin{eqnarray}\label{eqzrightsecond}
& &(x_{1}-z)^{l}(x_{2}-z)^{l}
\<Y^{R}_{Q(z)}(v,x_{2})Y^{R}_{Q(z)}(u,x_{1})\alpha,w\>\nonumber\\
&=&(x_{1}-z)^{l}(x_{2}-z)^{l}
\<\bar{Y}(u,x_{1})\bar{Y}(v,x_{2})\alpha,w\>.
\end{eqnarray}
Using (\ref{eqzrightfirst}), (\ref{eqzrightsecond}) and
the opposite Jacobi identity for $\bar{Y}$ we get
\begin{eqnarray}\label{eqzrightthird}
& &x_{0}^{-1}\delta\left(\frac{x_{1}-x_{2}}{x_{0}}\right)
(x_{1}-z)^{l}(x_{2}-z)^{l}
\<Y^{R}_{Q(z)}(u,x_{1})Y^{R}_{Q(z)}(v,x_{2})\alpha,w\>\nonumber\\
&-&x_{0}^{-1}\delta\left(\frac{x_{2}-x_{1}}{-x_{0}}\right)
(x_{1}-z)^{l}(x_{2}-z)^{l}
\<Y^{R}_{Q(z)}(v,x_{2})Y^{R}_{Q(z)}(u,x_{1})\alpha,w\>\nonumber\\
&=&x_{2}^{-1}\delta\left(\frac{x_{1}-x_{0}}{x_{2}}\right)
(x_{1}-z)^{l}(x_{2}-z)^{l}
\<\bar{Y}(Y(u,x_{0})v,x_{2})\alpha,w\>.
\end{eqnarray}

Let $m$ be any fixed integer. Since only finitely many
$u_{n}v$ are nonzero for $n\ge m$, there exists $l'\in {\N}$ 
(depending on $m$) such that
\begin{eqnarray}
(x_{2}-z)^{l'}\< Y^{R}_{Q(z)}(u_{n}v,x_{2})\alpha,w\>=
(x_{2}-z)^{l'}\< \bar{Y}(u_{n}v,x_{2})\alpha,w\>
\end{eqnarray}
for {\em all} $n\ge m$. Then
\begin{eqnarray}\label{eqzrightforth}
& &\Res_{x_{0}}x_{0}^{m}x_{2}^{-1}\delta\left(\frac{x_{1}-x_{0}}{x_{2}}\right)
(x_{2}-z)^{l'}\<\bar{Y}(Y(u,x_{0})v,x_{2})\alpha,w\>\nonumber\\
&=&\Res_{x_{0}}x_{0}^{m}x_{2}^{-1}\delta\left(\frac{x_{1}-x_{0}}{x_{2}}\right)
(x_{2}-z)^{l'}\< Y^{R}_{Q(z)}(Y(u,x_{0})v,x_{2})\alpha,w\>.
\end{eqnarray}
Let $l''=l+l'$. Combining (\ref{eqzrightthird}) with (\ref{eqzrightforth})
we get
\begin{eqnarray}\label{eqzrightfifth}
& &\Res_{x_{0}}x_{0}^{m}
x_{0}^{-1}\delta\left(\frac{x_{1}-x_{2}}{x_{0}}\right)
(x_{1}-z)^{l''}(x_{2}-z)^{l''}
\<Y^{R}_{Q(z)}(u,x_{1})Y^{R}_{Q(z)}(v,x_{2})\alpha,w\>\nonumber\\
&-&\Res_{x_{0}}x_{0}^{m}x_{0}^{-1}\delta\left(\frac{x_{2}-x_{1}}{-x_{0}}\right)
(x_{1}-z)^{l''}(x_{2}-z)^{l''}
\<Y^{R}_{Q(z)}(v,x_{2})Y^{R}_{Q(z)}(u,x_{1})\alpha,w\>\nonumber\\
&=&\Res_{x_{0}}x_{0}^{m}
x_{2}^{-1}\delta\left(\frac{x_{1}-x_{0}}{x_{2}}\right)
(x_{1}-z)^{l''}(x_{2}-z)^{l''}
\<Y^{R}_{Q(z)}(Y(u,x_{0})v,x_{2})\alpha,w\>.
\end{eqnarray}
Multiplying (\ref{eqzrightfifth}) by $(-z+x_{1})^{-l''}(-z+x_{2})^{-l''}$
we obtain
\begin{eqnarray}\label{eqzrightsixth}
& &\Res_{x_{0}}x_{0}^{m}
x_{0}^{-1}\delta\left(\frac{x_{1}-x_{2}}{x_{0}}\right)
\<Y^{R}_{Q(z)}(u,x_{1})Y^{R}_{Q(z)}(v,x_{2})\alpha,w\>\nonumber\\
&-&\Res_{x_{0}}x_{0}^{m}x_{0}^{-1}\delta\left(\frac{x_{2}-x_{1}}{-x_{0}}\right)
\<Y^{R}_{Q(z)}(v,x_{2})Y^{R}_{Q(z)}(u,x_{1})\alpha,w\>\nonumber\\
&=&\Res_{x_{0}}x_{0}^{m}
x_{2}^{-1}\delta\left(\frac{x_{1}-x_{0}}{x_{2}}\right)
\<Y^{R}_{Q(z)}(Y(u,x_{0})v,x_{2})\alpha,w\>.
\end{eqnarray}
Since $m$ was arbitrarily fixed, we may drop off $\Res_{x_{0}}x_{0}^{m}$
to obtain the Jacobi identity
for $Y^{R}_{Q(z)}$. $\;\;\;\;\Box$

In view of Proposition \ref{pqzleftrightcomm}, 
$Y^{L}_{Q(z)}$ and $Y^{R}_{Q(z)}$ give rise to a well defined
linear map 
\begin{eqnarray}
Y_{Q(z)}=Y^{L}_{Q(z)}\otimes Y^{R}_{Q(z)}:V\otimes V\rightarrow
(\End {\cal{D}}_{Q(z)}(W))[[x,x^{-1}]].
\end{eqnarray}
In particular,
\begin{eqnarray}
Y_{Q(z)}(u\otimes v,x)=Y_{Q(z)}^{L}(u,x)Y_{Q(z)}^{R}(v,x)
\end{eqnarray}
for $u,v\in V$. Combining Theorem \ref{tqzleft} with 
Propositions \ref{pqzright}
and \ref{pqzleftrightcomm}, using the arguments of Proposition 4.6.1 of [FHL]
 we immediately have:

\bt{tqzbimodule}
Let $W$  be a (left) weak $V$-module and let $z$ be a nonzero complex number. 
Then the pair $({\cal{D}}_{Q(z)}(W), Y_{Q(z)})$ carries 
the structure of a weak $V\otimes V$-module. $\;\;\;\;\Box$
\et

It is routine to check that the notion of ${\cal{D}}_{Q(z)}(W)$in the obvious way 
gives rise to a (contravariant) functor from the category of weak $V$-modules 
to the category of weak $V\otimes V$-modules.

\section{$V\otimes V$-homomorphisms and $Q(z)$-intertwining maps}
In this section, we shall prove that for $V$-modules 
$W,W_{1}$ and $W_{2}$, a $Q(z)$-intertwining map of type
${W'\choose W_{1}W_{2}}$ in the sense of [HL2] exactly amounts to a 
$V\otimes V$-homomorphism from $W_{1}\otimes W_{2}$ to
${\cal{D}}_{Q(z)}(W)$ and that a $Q(z)$-tensor product
of $V$-modules $W_{1}$ and $W_{2}$ in the sense of [HL2] 
exactly amounts to a universal from $W_{1}\otimes W_{2}$ to a
functor ${\cal{F}}_{Q(z)}$, which is essentially the composition of
the functor ${\cal{D}}_{Q(z)}$ with the contragredient module functor.

Following [HL2], for any $\C$-graded vector space 
$U=\coprod_{h\in {\C}}U_{(h)}$, we define the restricted dual
\begin{eqnarray}
U'=\coprod_{h\in {\C}}U_{(h)}^{*}
\end{eqnarray}
and the formal completion
\begin{eqnarray}
\overline{U}=\prod_{h\in \C}U_{(h)}.
\end{eqnarray}
Then
\begin{eqnarray}
\overline{U'}=U^{*}.
\end{eqnarray}
Furthermore, if $\dim U_{(h)}<\infty$ for $h\in \C$, e.g., $U=W$ is a $V$-module, 
we have
\begin{eqnarray}
\overline{U}=(U')^{*}.
\end{eqnarray}

Let $W_{1},W_{2}$ and $W_{3}$ be  $V$-modules.
A $Q(z)$-intertwining map of type ${W_{3}\choose W_{1}W_{2}}$ (see [HL2])
is a linear map $F: W_{1}\otimes W_{2}\rightarrow \overline{W_{3}}$ 
such that
\begin{eqnarray}\label{eqzjacobi}
& &z^{-1}\delta\left(\frac{x_{1}-x_{0}}{z}\right)
Y^{o}(v,x_{0})F(w_{(1)}\otimes w_{(2)})\nonumber\\
&=&x_{0}^{-1}\delta\left(\frac{x_{1}-z}{x_{0}}\right)
F(Y^{o}(v,x_{1})w_{(1)}\otimes w_{(2)})\nonumber\\
& &-x_{0}^{-1}\delta\left(\frac{z-x_{1}}{-x_{0}}\right)
F(w_{(1)}\otimes Y(v,x_{1})w_{(2)})
\end{eqnarray}
for $v\in V,\; w_{(1)}\in W_{1},\; w_{(2)}\in W_{2}$.
(With $W_{3}$ being a $V$-module,
the action of $Y^{o}(v,x)$ on $W_{3}$ extends to $\overline{W_{3}}$.)

Following [HL2], we choose $\log z$ so that
\begin{eqnarray}
\log z=\log |z|+i\arg z\;\mbox{ with }\; 0\le \arg z<2\pi.
\end{eqnarray}
Arbitrary values of the log function will be denoted by
\begin{eqnarray}
l_{p}(z)=\log z+2p\pi i
\end{eqnarray}
for $p\in \Z$.

We next recall from [HL2] a connection between 
intertwining operators of type ${W_{1}'\choose W_{3} W_{2}}$ 
and $Q(z)$-intertwining maps of type
${W_{3}'\choose W_{1}W_{2}}$.
Fix an integer $p$. Let ${\cal{Y}}$ be
any intertwining operator of type
${W_{1}'\choose W_{3}W_{2}}$. Define a linear map 
$F_{{\cal{Y}},p}$ from
$W_{1}\otimes W_{2}$ to $\overline{W_{3}'}=W_{3}^{*}$ by
\begin{eqnarray}
\<F_{{\cal{Y}},p}(w_{(1)}\otimes w_{(2)}),w_{(3)}\>_{W_{3}}
=\<{\cal{Y}}(w_{(3)},e^{l_{p}(z)})w_{(2)},w_{(1)}\>_{W_{1}}
\end{eqnarray}
for $w_{(1)}\in W_{1},\; w_{(2)}\in W_{2},\; w_{(3)}\in W_{3}$.
It was proved in [HL2] that $F_{{\cal{Y}},p}$ is a
$Q(z)$-intertwining map of type ${W_{3}'\choose W_{1}W_{2}}$.
Furthermore we have ([HL2], Proposition 4.7):

\bp{phlqzisom}
Let $W_{1},W_{2}$ and $W_{3}$ be (ordinary) $V$-modules. 
For $p\in {\Z}$, the correspondence 
${\cal{F}}_{p}: {\cal{Y}}\mapsto F_{{\cal{Y}},p}$
is a linear isomorphism from the space ${\cal{V}}_{W_{3}W_{2}}^{W_{1}'}$
of intertwining operators of the indicated type
to the space ${\cal{M}}[Q(z)]_{W_{1}W_{2}}^{W_{3}'}$ of 
$Q(z)$-intertwining maps of the indicated type.
\ep

Note that a $Q(z)$-intertwining map is a linear map from a
$V\otimes V$-module $W_{1}\otimes W_{2}$ to $\overline{W}$.
Let $A$ be a $V\otimes V$-module and $W$ a $V$-module. 
Then we may define a notion of $Q(z)$-intertwining map from $A$ to $\overline{W}$.
Denote by $Y_{1}$and $Y_{2}$ the two commuting (weak) $V$-module structures on $A$
obtained by identifying $V$ with $V\otimes {\C}$ and ${\C}\otimes V$, respectively.
A {\em $Q(z)$-intertwining map} from $A$ to $\overline{W}$
to be a linear map $F$ from $A$ to $\overline{W}$ such that
for $v\in V,\; a\in A$:
\begin{eqnarray}\label{eqzjacobiA}
& &z^{-1}\delta\left(\frac{x_{1}-x_{0}}{z}\right)
Y^{o}(v,x_{0})F(a)\nonumber\\
&=&z^{-1}\delta\left(\frac{x_{1}-x_{0}}{z}\right)
F(a)Y(v,x_{0})\nonumber\\
&=&x_{0}^{-1}\delta\left(\frac{x_{1}-z}{x_{0}}\right)
F(Y^{o}_{1}(v,x_{1})a)
-x_{0}^{-1}\delta\left(\frac{z-x_{1}}{-x_{0}}\right)
F(Y_{2}(v,x_{1})a).
\end{eqnarray}

The following is our key theorem of this section.

\bt{tqz}
Let $A$ be a $V\otimes V$-module, $W$ a $V$-module and 
$F$ a linear map from $A$ into $\overline{W}$ $(=(W')^{*})$. 
Then $F$ is a $Q(z)$-intertwining map from $A$ to $\overline{W}$ 
if and only if $F$ is a $V\otimes V$-homomorphism from $A$
into ${\cal{D}}_{Q(z)}(W')$ $(\subset (W')^{*})$.
\et

\pf Assume that $F$ is a $V\otimes V$-homomorphism from $A$
into ${\cal{D}}_{Q(z)}(W')$. For $v\in V,\; a\in A$, since
$F(a)\in {\cal{D}}_{Q(z)}(W')$, 
using (\ref{eqzrelation}) with $\alpha=F(a)$
and the fact that $F$ is a 
$V\otimes V$-homomorphism we obtain
\begin{eqnarray}
& &z^{-1}\delta\left(\frac{x_{1}-x_{0}}{z}\right)Y^{o}(v,x_{0})
F(a)\nonumber\\
&=&x_{0}^{-1}\delta\left(\frac{x_{1}-z}{x_{0}}\right)
(Y_{Q(z)}^{L})^{o}(v,x_{1})F(a)
-x_{0}^{-1}\delta\left(\frac{z-x_{1}}{-x_{0}}\right)
Y_{Q(z)}^{R}(v,x_{1})F(a)\nonumber\\
&=&x_{0}^{-1}\delta\left(\frac{x_{1}-z}{x_{0}}\right)
F(Y_{1}^{o}(v,x_{1})a)
-x_{0}^{-1}\delta\left(\frac{z-x_{1}}{-x_{0}}\right)
F(Y_{2}(v,x_{1})a).
\end{eqnarray}
This proves that $F$ is a $Q(z)$-intertwining map.

Conversely, assume that $F$ is a $Q(z)$-intertwining map.
Let $v\in V,\; a\in A$, and let
$l\in {\N}$ be such that $x^{l}Y_{2}(v,x)a\in A[[x]]$.
Taking $\Res_{x_{1}}x_{1}^{l}$ from (\ref{eqzjacobiA}) and using
the fundamental delta-function substitution properties we get
\begin{eqnarray}\label{eqzasso}
(x_{0}+z)^{l}Y^{o}(v,x_{0})F(a)
=(x_{0}+z)^{l}F(Y_{1}^{o}(v,x_{0}+z)a).
\end{eqnarray}
Since $Y_{1}^{o}(v,x_{0}+z)a\in A((x_{0}^{-1}))$, we have
\begin{eqnarray}
(x_{0}+z)^{l}F(Y_{1}^{o}(v,x_{0}+z)a)
\in \overline{W}((x_{0}^{-1})),
\end{eqnarray}
hence (from (\ref{eqzasso}))
\begin{eqnarray}
(x_{0}+z)^{l}Y^{o}(v,x_{0})F(a)
\in \overline{W}((x_{0}^{-1})).
\end{eqnarray}
In view of Lemma \ref{lqzequiv} we have
$F(a)\in {\cal{D}}_{Q(z)}(W')$.

Note that $F(a)Y(v,x_{0})w'=\<Y^{o}(v,x_{0})F(a),w'\>$ for $w'\in W'$.
Combining (\ref{eqzrelation}) with $\a=F(a)$ and (\ref{eqzjacobiA})
we obtain
\begin{eqnarray}\label{e4.13}
& &x_{0}^{-1}\delta\left(\frac{x_{1}-z}{x_{0}}\right)
(Y_{Q(z)}^{L})^{o}(v,x_{1})F(a)
-x_{0}^{-1}\delta\left(\frac{z-x_{1}}{-x_{0}}\right)
Y_{Q(z)}^{R}(v,x_{1})F(a)\nonumber\\
&=&x_{0}^{-1}\delta\left(\frac{x_{1}-z}{x_{0}}\right)
F(Y^{o}_{1}(v,x_{1})a)
-x_{0}^{-1}\delta\left(\frac{z-x_{1}}{-x_{0}}\right)
F(Y_{2}(v,x_{1})a).
\end{eqnarray}
Note that for $w'\in W'$,
\begin{eqnarray*}
& &\<(Y_{Q(z)}^{L})^{o}(v,x_{1})F(a),w'\>,\;\;
\<F(Y^{o}_{1}(v,x_{1})a),w'\>\in \C((x_{1}^{-1})),\\
& &\<Y_{Q(z)}^{R}(v,x_{1})F(a),w'\>,\;\;
\<F(Y_{2}(v,x_{1})a),w'\>\in \C((x)).
\end{eqnarray*}
Right after this theorem we shall prove a simple fact
(Lemma \ref{lsimplefactdelta}) 
which together with (\ref{e4.13}) immediately implies
\begin{eqnarray}
(Y_{Q(z)}^{L})^{o}(v,x_{0})F(a)&=&F(Y^{o}_{1}(v,x_{0})a)\\
Y_{Q(z)}^{R}(v,x_{1})F(a)&=&F(Y_{2}(v,x_{1})a).
\end{eqnarray}
This implies that $F$ is a $V\otimes V$-homomorphism.
$\;\;\;\;\Box$

We now prove the following simple fact which was used above:

\bl{lsimplefactdelta}
Let $z$ be a nonzero complex number and let 
$f(x)\in \C((x^{-1})),\; g(x)\in \C((x))$ be such that
\begin{eqnarray}\label{ef=g}
x_{0}^{-1}\delta\left(\frac{x_{1}-z}{x_{0}}\right)f(x_{1})-
x_{0}^{-1}\delta\left(\frac{z-x_{1}}{-x_{0}}\right)g(x_{1})=0.
\end{eqnarray}
Then $f(x)=g(x)=0$.
\el

\pf By taking $\Res_{x_{0}}$ from (\ref{ef=g}) we get $f(x_{1})=g(x_{1})$.
Consequently,
\begin{eqnarray}
f(x)=g(x)\in \C((x^{-1}))\cap \C((x))=\C[x,x^{-1}].
\end{eqnarray}
Then from (\ref{ef=g}) using the fundamental delta-function identities 
we get
\begin{eqnarray}
z^{-1}\delta\left(\frac{x_{1}-x_{0}}{z}\right)f(x_{1})=0.
\end{eqnarray}
Taking $\Res_{x_{1}}$ and using the fundamental delta-function properties 
we get $f(z+x_{0})=0.$ Thus $f(x)=0$ and $g(x)=f(x)=0$.
$\;\;\;\;\Box$

In view of  Proposition \ref{phlqzisom} and Theorem \ref{tqz} 
we immediately have:

\bc{cintertwininghomo}
Let $W, W_{1}$ and $W_{2}$ be (ordinary) $V$-modules and let $p\in \Z$. 
Then ${\cal{F}}_{p}$ gives rise to a linear isomorphism from
${\cal{V}}_{W'W_{2}}^{W_{1}'}$ to 
$\Hom_{V\otimes V}(W_{1}\otimes W_{2},{\cal{D}}_{Q(z)}(W'))$.
In particular,
\begin{eqnarray}
N_{W'W_{2}}^{W_{1}'}=\dim {\cal{V}}_{W'W_{2}}^{W_{1}'}
=\dim \Hom_{V\otimes V}(W_{1}\otimes W_{2},{\cal{D}}_{Q(z)}(W')).\;\;\;\;\Box
\end{eqnarray}
\ec

Recall a notion from [HL2].
A {\em generalized} $V$-module is a weak $V$-module on which
$L(0)$ acts semisimply. 
Then for a generalized $V$-module $W$ we have 
the $L(0)$-eigenspace decomposition: $W=\coprod_{h\in {\C}}W_{(h)}$.
Thus, a generalized $V$-module 
satisfies all the axioms defining the notion of a $V$-module
([FLM], [FHL]) except the two grading restrictions on the 
homogeneous subspaces.

Let $W$ be a $V$-module. We denote by $R_{Q(z)}(W)$ the sum of 
(ordinary) $V\otimes V$-submodules of ${\cal{D}}_{Q(z)}(W)$, on which
$L(0)\otimes 1$ acts semisimply. Then $R_{Q(z)}(W)$ is a generalized
$V\otimes V$-module. It is clear that Theorem \ref{tqz} and 
Corollary \ref{cintertwininghomo} still hold 
when ${\cal{D}}_{Q(z)}(W')$ is replaced by $R_{Q(z)}(W')$.

Note that the functor
${\cal{F}}^{c}$ with ${\cal{F}}^{c}(W)=W'$ is 
a contravariant functor from the category of $V$-modules to itself
and that ${\cal{D}}_{Q(z)}$ ($R_{Q(z)}$) 
is a contravariant functor from the category of $V$-modules
to the category of weak (generalized) $V\otimes V$-modules.
The composition of the two (contravariant) functors 
gives us a (covariant) functor from the category of 
$V$-modules to the category of weak (generalized) $V\otimes V$-modules.

\bd{dffunctor}
{\em Define ${\cal{F}}_{Q(z)}$ to be the 
(covariant) functor from the category of 
$V$-modules to the category of generalized $V\otimes V$-modules
such that
\begin{eqnarray}
& &{\cal{F}}_{Q(z)}(W)=R_{Q(z)}(W')\;\;\;\mbox{ for }
W\in \ob\; {\cal{C}}_{V},\\
& &{\cal{F}}_{Q(z)}(f)=(f')^{*}|_{R_{Q(z)}(W_{1}')}
=\bar{f}|_{R_{Q(z)}(W_{1}')}\;\;\;\mbox{ for }
f\in \Hom_{V}(W_{1},W_{2}),\; W_{1},W_{2}\in \ob\; {\cal{C}}_{V},
\end{eqnarray}
where ${\cal{C}}_{V}$ by definition denotes the category of $V$-modules.}
\ed

Next, we shall give a connection between the functor ${\cal{D}}_{Q(z)}$
and the notion of the $Q(z)$-tensor product defined in [HL2-3].
Recall from [HL2] that
a {\em $Q(z)$-product of $W_{1}$ and $W_{2}$} is a 
$V$-module $W_{3}$
together with a $Q(z)$-intertwining map $F$ of type 
${W_{3}\choose W_{1}W_{2}}$ and it is denoted by $(W_{3},F)$.
Let $(W_{3},F)$ and $(W_{4},G)$ be two $Q(z)$-products of 
$W_{1}$ and $W_{2}$. A {\em morphism} from 
$(W_{3},F)$ to $(W_{4},G)$ is a $V$-module map $\eta$ from $W_{3}$
to $W_{4}$ such that
\begin{eqnarray}
G=\bar{\eta}\circ F
\end{eqnarray}
where $\bar{\eta}$ is the natural map from $\overline{W}_{3}$ to 
$\overline{W}_{4}$ uniquely extending $\eta$.

\bd{dhlpzpro}
{\em [HL2,4] A {\em $Q(z)$-tensor product of $W_{1}$ and $W_{2}$} 
is a $Q(z)$-product
$(W,F)$ such that for any $Q(z)$-product $(W_{3},G)$, there is a unique
morphism from $(W,F)$ to $(W_{3},G)$. The $V$-module $W$ is called
a {\em $Q(z)$-tensor product module} of $W_{1}$ and $W_{2}$.}
\ed

In view of Theorem \ref{tqz}, for a $Q(z)$-product $(W,F)$ 
of $V$-modules $W_{1}$ and $W_{2}$,
$F$ is a morphism from $W_{1}\otimes W_{2}$
to $R_{Q(z)}(W')$ in the category
of generalized $V\otimes V$-modules. Then the condition for a $Q(z)$-product
$(W,F)$ of $W_{1}$ and $W_{2}$ being a $Q(z)$-tensor product
amounts to that for any $V$-module $W_{3}$ and any 
$V\otimes V$-homomorphism $G$ from $W_{1}\otimes W_{2}$ to
$R_{Q(z)}(W_{3}')={\cal{F}}_{Q(z)}(W_{3})$ there exists a unique
$V$-homomorphism $\eta$ from $W$ to $W_{3}$ such that
\begin{eqnarray}
G=\bar{\eta}\circ F=(\eta')^{*}\circ F=
{\cal{F}}_{Q(z)}(\eta)\circ F.
\end{eqnarray}
The later condition, in terms of a categorical notion (cf. [J]),
exactly amounts to that $(W,F)$ is 
a universal from $W_{1}\otimes W_{2}$ to the functor
${\cal{F}}_{Q(z)}$.

To summarize we have:

\bp{puniversalobjects}
Let $W_{1}, W_{2}$ be $V$-modules and let $(W,F)$ be 
a $Q(z)$-product of $W_{1}$ and $W_{2}$.
Then $(W,F)$ is a $Q(z)$-tensor product of $W_{1}$ and $W_{2}$
if and only if 
$(W,F)$ is a universal from $W_{1}\otimes W_{2}$ to the functor
${\cal{F}}_{Q(z)}$.
$\;\;\;\;\Box$
\ep

It was proved in [Li2] that if every $V$-module is completely
reducible, then every $V\otimes V$-module on which $L(0)\otimes 1$
acts semisimply is completely reducible, in particular,
$R_{Q(z)}(W)$ is completely reducible for every $V$-module $W$.
Let $S$ be a fixed complete set of representatives of equivalent classes
of irreducible $V$-modules. Then it follows from [FHL] that
$\{ W_{1}\otimes W_{2}\;|\; W_{1},W_{2}\in S\}$ is 
a complete set of representatives of equivalent classes
of irreducible $V\otimes V$-modules.
Combining this with Corollary \ref{cintertwininghomo}, 
we immediately have:

\bp{prqzwmodule}
Assume that $V$ is rational in the sense of Huang-Lepowsky. 
Let $W$ be a $V$-module and $z$ a nonzero complex number. 
Then $R_{Q(z)}(W)$ is a direct sum of finitely many irreducible
$V\otimes V$-modules. In particular, $R_{Q(z)}(W)$ is 
an (ordinary) $V\otimes V$-module.$\;\;\;\;\Box$
\ep

Let $S$ be a fixed complete set of representatives of equivalent classes
of irreducible $V$-modules. For a $V$-module $W$, we set
\begin{eqnarray}
H_{Q(z)}(S,W)=\coprod_{W_{1},W_{2}\in S}{\cal{M}}[Q(z)]_{W_{1}W_{2}}^{W}
\otimes (W_{1}\otimes W_{2}).
\end{eqnarray}
We consider $H_{Q(z)}(S,W)$ as a generalized $V\otimes V$-module in the obvious way.
In view of Theorem \ref{tqz} we obtain a $V\otimes V$-homomorphism
\begin{eqnarray}
\Psi_{W}[Q(z)]: H_{Q(z)}(S,W)\rightarrow R_{Q(z)}(W')
\end{eqnarray}
such that
\begin{eqnarray}
\Psi_{W}[Q(z)](F\otimes w_{(1)}\otimes w_{(2)})
=F(w_{(1)}\otimes w_{(2)})
\end{eqnarray}
for $F\in {\cal{M}}[Q(z)]_{W_{1}W_{2}}^{W},\; w_{(1)}\in W_{1},
\; w_{(2)}\in W_{2}$ with $W_{1},W_{2}\in S$.
Then the same argument of Theorem 4.17 in [Li2] gives:

\bp{pQ(z)decomposition}
Assume that every $V$-module is completely reducible.
Then for any $V$-module $W$, $\Psi_{W}[Q(z)]$ is a $V\otimes V$-isomorphism
onto $R_{Q(z)}(W')$. In particular, if $V$ is rational in the sense of 
Huang-Lepowsky, then this is true and $R_{Q(z)}(W')$ is a $V\otimes V$-module.
$\;\;\;\;\Box$
\ep

Set
\begin{eqnarray}
H(S)=\coprod_{W\in S}W'\otimes W.
\end{eqnarray}
Since $Y$ is an intertwining operator of type ${W\choose VW}$
for any $V$-module $W$, by Proposition \ref{phlqzisom} 
([HL2], Proposition 4.7) the linear map $F_{W}$ defined by
\begin{eqnarray}
\<F_{W}(w'\otimes w),v\>=\<w', Y(v,z)w\>
\end{eqnarray}
for $w'\in W',\; w\in W,\;v\in V,\;v\in V$ is a $Q(z)$-intertwining map
of type ${V'\choose W' W}$. It follows from Theorem \ref{tqz} that $F_{W}$
is a $V\otimes V$-homomorphism from $W'\otimes W$
to $R_{Q(z)}(V)$.
Then we have a $V\otimes V$-homomorphism $\Psi[Q(z)]$ from $H(S)$
to $R_{Q(z)}(V)$ such that
\begin{eqnarray}
\<\Psi[Q(z)](w'\otimes w),v\>=\<w', Y(v,z)w\>
\end{eqnarray}
for $v\in V,\; w'\in W',\; w\in W$ with $W\in S$.
It is easy to show that ${\cal{V}}_{VW}^{W}$ is
one-dimensional for any irreducible $V$-module $W$.
Then we immediately have:

\bp{pregularQ(z)decomposition}
Assume that every $V$-module is completely reducible.
Then $R_{Q(z)}(V)$ is a $V\otimes V$-module
and $\Psi[Q(z)]$ is a $V\otimes V$-isomorphism
onto $R_{Q(z)}(V)$. In particular, the assertion holds
if $V$ is rational in the sense of Huang-Lepowsky.$\;\;\;\;\Box$
\ep

\br{rfusionrule}
{\em We here give a detailed proof for ${\cal{V}}_{VW}^{W}$ being
one-dimensional for any irreducible $V$-module $W$.
First, we prove that ${\cal{V}}_{VW_{1}}^{W_{2}}$ is linearly
isomorphic to
$\Hom_{V}(W_{1},W_{2})$ for any $V$-modules
$W_{1}$ and $W_{2}$. Then it follows from Schur lemma (see [FHL])
that ${\cal{V}}_{VW}^{W}$ is one-dimensional for any 
irreducible $V$-module $W$.
Let
${\cal{Y}}$ be an intertwining operator of type ${W_{2}\choose VW_{1}}$.
Then ${\cal{Y}}({\bf 1},x)\in \Hom (W_{1},W_{2})$ because
$${d\over dx}{\cal{Y}}({\bf 1},x)={\cal{Y}}(L(-1){\bf 1},x)=0.$$
Since $v_{i}{\bf 1}=0$ for $v\in V,\; i\ge 0$, from the Jacobi identity
we get
$$Y(v,x_{1}){\cal{Y}}({\bf 1},x)={\cal{Y}}({\bf 1},x)Y(v,x_{1}).$$
That is, ${\cal{Y}}({\bf 1},x)\in \Hom_{V}(W_{1},W_{2})$. Then we have a
linear map $E: {\cal{Y}}\mapsto {\cal{Y}}({\bf 1},x)$ from 
${\cal{V}}_{VW_{1}}^{W_{2}}$ to $\Hom_{V}(W_{1},W_{2})$.
If ${\cal{Y}}({\bf 1},x)=0$, then
it follows from the Jacobi identity of ${\cal{Y}}$ and and creation property
$v_{-1}{\bf 1}=v$  that ${\cal{Y}}(v,x)=0$ for all $v\in V$.
Thus, $E$ is injective.
On the other hand, let $f$ be a $V$-homomorphism from $W_{1}$ to $W_{2}$.
Then it is clear that $f\circ Y$ is an intertwining operator of 
type ${W_{2}\choose VW_{1}}$ such that 
$$E(f\circ Y)=(f\circ Y)({\bf 1},x)=fY({\bf 1},x)=f.$$
Then $E$ is onto. Therefore, $E$ is a linear isomorphism.}
\er

\br{rpzqzconnection}
{\em Let $W_{1}, W_{2}$ and $W_{3}$ be (ordinary) $V$-modules.
{}From the proof of Propositions 4.7 and 12.2 in [HL2,4],
$Q(z)$-intertwining maps
of type ${W_{3}\choose W_{1}W_{2}}$
and $P(z)$-intertwining
maps of type ${W_{1}'\choose W_{3}'W_{2}}$
can be canonically related as follows:
Let $F$ be a linear map from $W_{1}\otimes W_{2}$ 
to $\overline{W_{3}}$. Define a 
linear map $F'$ from $W_{3}'\otimes W_{2}$ 
to $\overline{W_{1}'}=W_{1}^{*}$ by
\begin{eqnarray}
\<F'(w_{(3)}'\otimes w_{(2)}),w_{(1)}\>_{W_{1}}
=\<w_{(3)}',F(w_{(1)}\otimes w_{(2)})\>_{W_{3}}
\end{eqnarray}
for $w_{(1)}\in W_{1},\; w_{(2)}\in W_{2},\; w_{(3)}'\in W_{3}'$.
Then it is straightforward to check that
$F$ is a $Q(z)$-intertwining map of type ${W_{3}\choose W_{1}W_{2}}$
if and only if $F'$ is a $P(z)$-intertwining map of type 
${W_{1}'\choose W_{3}'W_{2}}$.}
\er


\section{Huang-Lepowsky's $Q(z)$-tensor functor}

In this section, we shall extend Huang-Lepowsky's $Q(z)$-tensor functor
constructed in [HL2-3] for the category of $V$-modules 
to a functor $T_{Q(z)}$ from the category of 
$V\otimes V$-modules to the category of $V$-modules.
We then show that the functor ${\cal{F}}_{Q(z)}$ constructed in
Section 4 is a right adjoint of $T_{Q(z)}$.
As mentioned in Introduction, in establishing the functor 
$T_{Q(z)}$ one may follow [HL2-3] and presumably
the main result (Theorem \ref{texthl}) follows from 
Huang-Lepowsky's arguments in [HL2-3] without any
further efforts. We here choose to use a different approach by 
making use of the regular representations.

Let $A$ be a weak $V\otimes V$-module for now and
let $A^{*}$ be the algebraic dual of $A$.
We define a vertex operator map $Y_{Q(z)}'$ from $V$
to $(\End A^{*})[[x,x^{-1}]]$ by (cf. [HL2], (5.4))
\begin{eqnarray}\label{eqzaction}
& &\<Y_{Q(z)}'(v,x_{0})\lambda, a\>\nonumber\\
&=&\Res_{x_{1}}x_{0}^{-1}\delta\left(\frac{x_{1}-z}{x_{0}}\right)
\<\lambda,Y_{1}^{o}(v,x_{1})a\>
-\Res_{x_{1}}x_{0}^{-1}\delta\left(\frac{z-x_{1}}{-x_{0}}\right)
\<\lambda,Y_{2}(v,x_{1})a\>\hspace{1cm}
\end{eqnarray}
for $v\in V,\; \lambda\in A^{*},\; a\in A$.

\bd{dQ(z)compatibility} 
{\em An element $\lambda$ of $A^{*}$ is said to satisfy
Huang-Lepowsky's {\em $Q(z)$-compatibility condition} (see [HL2-3]) if

(a) the lower truncation condition holds for all $v\in V$:
\begin{eqnarray}\label{eqzlower}
Y_{Q(z)}'(v,x)\lambda\in A^{*}((x)).
\end{eqnarray}

(b) the following Jacobi identity relation holds for all $v\in V,\;a\in A$:
\begin{eqnarray}\label{eqzcompatibility}
& &z^{-1}\delta\left(\frac{x_{1}-x_{0}}{z}\right)
\<Y_{Q(z)}'(v,x_{0})\lambda, a\>\nonumber\\
&=&x_{0}^{-1}\delta\left(\frac{x_{1}-z}{x_{0}}\right)
\<\lambda,Y_{1}^{o}(v,x_{1})a\>
-x_{0}^{-1}\delta\left(\frac{z-x_{1}}{-x_{0}}\right)
\<\lambda,Y_{2}(v,x_{1})a\>.\hspace{1cm}
\end{eqnarray}}
\ed

Clearly, all the linear functionals on $A$ that satisfy 
the Huang-Lepowsky's $Q(z)$-compatibility condition form a subspace 
of $A^{*}$, which we denote by $\tilde{T}_{Q(z)}'(A)$.

\br{rpresume}
{\em Note that if $A=W_{1}\otimes W_{2}$ for $V$-modules $W_{1},W_{2}$, then
$Y'_{Q(z)}(v,x_{0})$ is exactly the one defined in [HL2] and 
the Huang-Lepowsky's $Q(z)$-compatibility condition is exactly the 
one defined therein.}
\er

Let $v\in V,\; \lambda\in A^{*},\; a\in A$ and let
$l\in \N$ be such that $x_{1}^{l}Y_{2}(v,x_{1})a\in A[[x_{1}]]$.
Then from (\ref{eqzaction}) we get (cf. [DL])
\begin{eqnarray}\label{e2.53}
(x_{0}+z)^{l}\<Y'_{Q(z)}(v,x_{0})\lambda,a\>
=(x_{0}+z)^{l}\<\lambda,Y_{1}^{o}(v,x_{0}+z)a\>.
\end{eqnarray}
Since the expression on the right-hand side lies in $\C((x_{0}^{-1}))$,
so does the expression on the left-hand side.
Furthermore, if $\lambda\in \tilde{T}'_{Q(z)}(A)$,
using the lower truncation condition we obtain
\begin{eqnarray}
(x_{0}+z)^{l}\<Y'_{Q(z)}(v,x_{0})\lambda,a\>\in 
\C((x_{0}))\cap \C((x_{0}^{-1}))=\C[x_{0},x_{0}^{-1}].
\end{eqnarray}
Consequently,
the formal series $\<Y'_{Q(z)}(v,x_{0})\lambda,a\>$
converges to a rational function in the subring 
$\C[x_{0}, x_{0}^{-1}, (x_{0}+z)^{-1}]$.

The following is an equivalent form of the Huang-Lepowsky's
$Q(z)$-compatibility condition:

\bl{lweakcommqzcomp}
An element $\lambda$ of $A^{*}$ satisfies the Huang-Lepowsky's
$Q(z)$-compatibility condition if and only if 
for every $v\in V$, there exists $k\in \N$ such that
\begin{eqnarray}\label{eqzweakcomm}
(x_{1}-z)^{k}\<\lambda,Y_{1}^{o}(v,x_{1})a\>
=(x_{1}-z)^{k}\<\lambda,Y_{2}(v,x_{1})a\>
\end{eqnarray}
for all $a\in A$.
\el

\pf Assume that $\lambda$ satisfies the Huang-Lepowsky's $Q(z)$-compatibility 
condition. For $v\in V$, let $k$ be a nonnegative integer such that
$x^{k}Y_{Q(z)}'(v,x)\lambda\in A^{*}[[x]]$.
By applying $\Res_{x_{0}}x_{0}^{k}$ to (\ref{eqzcompatibility})
we obtain (\ref{eqzweakcomm}).

Conversely, assume (\ref{eqzweakcomm}).
It follows from (\ref{eqzaction}) and (\ref{eqzweakcomm}) that
$$x_{0}^{k}Y_{Q(z)}'(v,x_{0})\lambda\in A^{*}[[x_{0}]].$$
That is, $\lambda$ satisfies the lower truncation condition.
Furthermore, 
it is well known (cf. [DL, Li1]) that
(\ref{eqzweakcomm}) and (\ref{eqzaction}) imply
(\ref{eqzcompatibility}). This proves that $\lambda$
satisfies the Huang-Lepowsky's $Q(z)$-compatibility 
condition.$\;\;\;\;\Box$

Note that if (\ref{eqzweakcomm}) holds for some $k$, then (\ref{eqzweakcomm}) 
holds with $k$ being replaced by any integer greater than $k$.

\br{rrationalityhlQ(z)}
{\em With Lemma \ref{lweakcommqzcomp}, it follows from
[FLM], [FHL], [DL] that the formal series
$\<\lambda, Y_{1}^{o}(v,x)a\>$ and 
$\<\lambda, Y_{2}(v,x)a\>$ absolutely converge 
in the domains $|x|>|z|$ and $0<|x|<|z|$, respectively,
to a common rational function with only two possible poles at
$x=0, z$ such that the order of the pole at $x=z$ is universally bounded
for all $a\in A$. }
\er

\br{rhlQ(z)Drelation}
{\em Note that the order of the pole at $x=0$ for the rational function
discussed in Remark \ref{rrationalityhlQ(z)} in general depends on
$a\in A$. So, an element $\lambda$ of $\tilde{T}_{Q(z)}'(A)$ is not quite
an element of ${\cal{D}}_{P(z)}(A,Y_{1})$, nor an element of 
${\cal{D}}_{Q(-z)}(A,Y_{2})$. On the other hand,
as illustrated in the proof of Theorem \ref{texthl}, $\lambda$
can be considered as an element of ${\cal{D}}_{P(z)}(U,Y_{1})$
for any finitely generated $V$-submodule $U$ of $(A,Y_{1})$. }
\er

We have (cf. [HL2], Proposition 6.2):

\bp{ptqzstable}
The subspace  $\tilde{T}_{Q(z)}'(A)$ of $A^{*}$ is stable
under the action of all components of $Y_{Q(z)}'(v,x_{0})$ for $v\in V$, i.e., 
\begin{eqnarray}
Y_{Q(z)}'(v,x_{0})\lambda\in \tilde{T}_{Q(z)}'(A)[[x_{0},x_{0}^{-1}]]
\;\;\;\mbox{ for }\lambda\in \tilde{T}_{Q(z)}'(A).
\end{eqnarray}
\ep

\pf Let $\lambda\in \tilde{T}_{Q(z)}'(A)$ and let
$u,v\in V$. Let $k\in \N$ be such that all the following identities hold
for all $a\in A$:
\begin{eqnarray}
& &(x_{1}-y)^{k}Y_{1}^{o}(u,x_{1})Y_{1}^{o}(v,y)a
=(x_{1}-y)^{k}Y_{1}^{o}(v,y)Y_{1}^{o}(u,x_{1})a,\\
& &(x_{1}-y)^{k}Y_{2}(u,x_{1})Y_{2}(v,y)a
=(x_{1}-y)^{k}Y_{2}(v,y)Y_{2}(u,x_{1})a,\\
& &(y-z)^{k}\<\lambda,Y_{1}^{o}(v,y)a\>=
(y-z)^{k}\<\lambda,Y_{2}(v,y)a\>,\label{e5.10}
\end{eqnarray}
which follow from weak commutativity and Lemma \ref{lweakcommqzcomp}.
Replacing $a$ with $Y_{1}^{o}(u,x_{1})a$ and $Y_{2}(u,x_{1})a$ in (\ref{e5.10})
we get
\begin{eqnarray}
& &(y-z)^{k}\<\lambda,Y_{1}^{o}(v,y)Y_{1}^{o}(u,x_{1})a\>=
(y-z)^{k}\<\lambda,Y_{2}(v,y)Y_{1}^{o}(u,x_{1})a\>,\\
& &(y-z)^{k}\<\lambda,Y_{1}^{o}(v,y)Y_{2}(u,x_{1})a\>=
(y-z)^{k}\<\lambda,Y_{2}(v,y)Y_{2}(u,x_{1})a\>.
\end{eqnarray}
Noting that the actions $Y_{1}$ and $Y_{2}$ of $V$ commute,
using the fundamental properties of delta function and
these identities we get
\begin{eqnarray}\label{e3.44}
& &(x_{0}+z-y)^{k}(y-z)^{k}
\<Y_{Q(z)}'(u,x_{0})\lambda, Y_{1}^{o}(v,y)a\>\nonumber\\
&=&\Res_{x_{1}}x_{0}^{-1}\delta\left(\frac{x_{1}-z}{x_{0}}\right)
(x_{1}-y)^{k}(y-z)^{k}
\<\lambda,Y_{1}^{o}(u,x_{1})Y_{1}^{o}(v,y)a\>\nonumber\\
& &-\Res_{x_{1}}x_{0}^{-1}\delta\left(\frac{z-x_{1}}{-x_{0}}\right)
(x_{1}-y)^{k}(y-z)^{k}\<\lambda,Y_{2}(u,x_{1})Y_{1}^{o}(v,y)a\>\nonumber\\
&=&\Res_{x_{1}}x_{0}^{-1}\delta\left(\frac{x_{1}-z}{x_{0}}\right)
(x_{1}-y)^{k}(y-z)^{k}
\<\lambda,Y_{1}^{o}(v,y)Y_{1}^{o}(u,x_{1})a\>\nonumber\\
& &-\Res_{x_{1}}x_{0}^{-1}\delta\left(\frac{z-x_{1}}{-x_{0}}\right)
(x_{1}-y)^{k}(y-z)^{k}\<\lambda,Y_{1}^{o}(v,y)Y_{2}(u,x_{1})a\>\nonumber\\
&=&\Res_{x_{1}}x_{0}^{-1}\delta\left(\frac{x_{1}-z}{x_{0}}\right)
(x_{1}-y)^{k}(y-z)^{k}
\<\lambda,Y_{2}(v,y)Y_{1}^{o}(u,x_{1})a\>\nonumber\\
& &-\Res_{x_{1}}x_{0}^{-1}\delta\left(\frac{z-x_{1}}{-x_{0}}\right)
(x_{1}-y)^{k}(y-z)^{k}\<\lambda,Y_{2}(v,y)Y_{2}(u,x_{1})a\>\nonumber\\
&=&\Res_{x_{1}}x_{0}^{-1}\delta\left(\frac{x_{1}-z}{x_{0}}\right)
(x_{1}-y)^{k}(y-z)^{k}
\<\lambda,Y_{1}^{o}(u,x_{1})Y_{2}(v,y)a\>\nonumber\\
& &-\Res_{x_{1}}x_{0}^{-1}\delta\left(\frac{z-x_{1}}{-x_{0}}\right)
(x_{1}-y)^{k}(y-z)^{k}\<\lambda,Y_{2}(u,x_{1})Y_{2}(v,y)a\>\nonumber\\
&=&(x_{0}+z-y)^{k}(y-z)^{k}\<Y_{Q(z)}'(u,x_{0})\lambda,Y_{2}(v,y)a\>.
\end{eqnarray}
Let $n\in \Z$ be arbitrarily fixed. Let $r\in \N$ be such that
$$x_{0}^{n+r}Y_{Q(z)}'(u,x_{0})\lambda\in A^{*}[[x_{0}]].$$
Then using (\ref{e3.44}) we get
\begin{eqnarray}
& &(y-z)^{r+2k}\Res_{x_{0}}x_{0}^{n}
\<Y_{Q(z)}'(u,x_{0})\lambda, Y_{1}^{o}(v,y)a\>\nonumber\\
&=&\sum_{i=0}^{r+k}{r+k\choose i}\Res_{x_{0}}(y-z-x_{0})^{r+k-i}x_{0}^{n+i}
(y-z)^{k}\<Y_{Q(z)}'(u,x_{0})\lambda, Y_{1}^{o}(v,y)a\>\nonumber\\
&=&\sum_{i=0}^{r}{r+k\choose i}\Res_{x_{0}}(y-z-x_{0})^{r+k-i}x_{0}^{n+i}
(y-z)^{k}\<Y_{Q(z)}'(u,x_{0})\lambda, Y_{1}^{o}(v,y)a\>\nonumber\\
&=&\sum_{i=0}^{r}{r+k\choose i}\Res_{x_{0}}(y-z-x_{0})^{r+k-i}x_{0}^{n+i}
(y-z)^{k}\<Y_{Q(z)}'(u,x_{0})\lambda,Y_{2}(v,y)a\>\nonumber\\
&=&\Res_{x_{0}}\sum_{i=0}^{r+k}{r+k\choose i}(y-z-x_{0})^{r+k-i}x_{0}^{n+i}
(y-z)^{k}\<Y_{Q(z)}'(u,x_{0})\lambda,Y_{2}(v,y)a\>\nonumber\\
&=&(y-z)^{r+2k} \Res_{x_{0}}x_{0}^{n}
\<Y_{Q(z)}'(u,x_{0})\lambda,Y_{2}(v,y)a\>.
\end{eqnarray}
Then from Lemma \ref{lweakcommqzcomp} we get
\begin{eqnarray}
\Res_{x_{0}}x_{0}^{n}Y_{Q(z)}'(u,x_{0})\lambda\in 
\tilde{T}_{Q(z)}'(A).
\end{eqnarray}
This completes the proof. $\;\;\;\;\Box$

We have (cf. [HL2], Theorem 6.1):

\bt{texthl}
The pair $(\tilde{T}_{Q(z)}'(A),Y_{Q(z)}')$ carries the structure
of a weak $V$-module.
\et

\pf From the $Q(z)$-compatibility condition and 
Proposition \ref{ptqzstable} we have
$$Y_{Q(z)}'(v,x)\lambda\in \tilde{T}_{Q(z)}'(A)((x))
\;\;\;\mbox{ for }v\in V,\; \lambda \in \tilde{T}_{Q(z)}'(A).$$
{}From (\ref{eqzaction}), clearly, $Y_{Q(z)}'({\bf 1},x)=1$. 
Then it remains to prove the Jacobi identity. 

Let $a\in A$. Let
$U$ be the $V$-submodule of $(A,Y_{1})$ generated by $a$.
Denote by $p_{U}$ the natural projection from $A^{*}$ to $U^{*}$.
For $v\in V$, let $l\in \N$ be such that $x^{l}Y_{2}(v,x)a\in A[[x]]$.
Since $U$ is generated from $a$ by $Y_{1}$ and since $Y_{1}$ and $Y_{2}$ 
commute, we have
\begin{eqnarray}
x^{l}Y_{2}(v,x)b\in A[[x]]\;\;\;\mbox{ for all }b\in U.
\end{eqnarray}
For $\lambda\in \tilde{T}_{Q(z)}'(A)$, let $k\in \N$ be such that 
 (\ref{eqzweakcomm}) holds.
Then 
\begin{eqnarray}
x_{1}^{l}(x_{1}-z)^{k}\<\lambda,Y_{1}^{o}(v,x_{1})b\>
=x_{1}^{l}(x_{1}-z)^{k}\<\lambda,Y_{2}(v,x_{1})b\>
\in \C[x_{1}]
\end{eqnarray}
for all $b\in U$. Thus $p_{U}(\lambda)\in {\cal{D}}_{P(z)}(U,Y_{1})$.
Recall (\ref{e2.53}) and (\ref{e2.26}):
\begin{eqnarray}
& &(x_{0}+z)^{l}\<Y_{Q(z)}'(v,x_{0})\lambda,b\>
=(x_{0}+z)^{l}\<\lambda,Y_{1}^{o}(v,x_{0}+z)b\>,\\
& &(x_{0}+z)^{l}\<Y_{P(z)}^{L}(v,x_{0})p_{U}(\lambda),b\>
=(x_{0}+z)^{l}\<\lambda,Y_{1}^{o}(v,x_{0}+z)b\>.
\end{eqnarray}
Then
\begin{eqnarray}
(x_{0}+z)^{l}\<Y_{Q(z)}'(v,x_{0})\lambda,b\>
=(x_{0}+z)^{l}\<Y_{P(z)}^{L}(v,x_{0})p_{U}(\lambda),b\>,
\end{eqnarray}
which by multiplying both sides by $(z+x_{0})^{-l}$ gives
\begin{eqnarray}
\<Y_{Q(z)}'(v,x_{0})\lambda,b\>=\<Y_{P(z)}^{L}(v,x_{0})p_{U}(\lambda),b\>.
\end{eqnarray}
That is,
\begin{eqnarray}
p_{U}(Y_{Q(z)}'(v,x_{0})\lambda)=Y_{P(z)}^{L}(v,x_{0})p_{U}(\lambda)
\end{eqnarray}
for $v\in V,\; \lambda \in \tilde{T}_{Q(z)}(A)$.
Then for $u,v\in V,\;b\in U$,
\begin{eqnarray}
\<Y_{Q(z)}'(u,x_{1})Y_{Q(z)}'(v,x_{2})\lambda,b\>
&=&\<p_{U}(Y_{Q(z)}'(u,x_{1})Y_{Q(z)}'(v,x_{2})\lambda),b\>\nonumber\\
&=&\<Y_{P(z)}^{L}(u,x_{1})Y^{L}_{P(z)}(v,x_{2})p_{U}(\lambda),b\>.
\end{eqnarray}
(We are using Proposition \ref{ptqzstable}.) Similarly,
\begin{eqnarray}
& &\<Y_{Q(z)}'(v,x_{2})Y_{Q(z)}'(u,x_{1})\lambda,b\>
=\<Y^{L}_{P(z)}(v,x_{2})Y_{P(z)}^{L}(u,x_{1})p_{U}(\lambda),b\>,\\
& &\<Y_{Q(z)}'(Y(u,x_{0})v,x_{2})\lambda,b\>
=\<Y_{P(z)}^{L}(Y(u,x_{0})v,x_{2})p_{U}(\lambda),b\>.
\end{eqnarray}
Then the Jacobi identity 
for $Y_{Q(z)}'$ applied to any element of $U$, in particular, to $a$,
follows from the Jacobi identity for $Y_{P(z)}^{L}$.
$\;\;\;\;\Box$

Denote by $T_{Q(z)}'(A)$ the sum of all (ordinary) $V$-submodules of 
$\tilde{T}_{Q(z)}'(A)$. Then $T_{Q(z)}'(A)$ is a generalized $V$-module.

\br{rclassicalfact}
{\em Let $U_{1},U_{2}$ be arbitrary vector spaces.
Let $f$ be a linear map from $U_{1}$ to $U_{2}^{*}$. Then
$f^{*}$ is a linear map from $(U_{2}^{*})^{*}$ to $U_{1}^{*}$, hence
$f^{*}|_{U_{2}}$ is a linear map from $U_{2}$ to $U_{1}^{*}$.
A simple fact is that the map $f\mapsto f^{*}|_{U_{2}}$ is a 
linear isomorphism from $\Hom (U_{1},U_{2}^{*})$  to
$\Hom (U_{2},U_{1}^{*})$ with
\begin{eqnarray}
(f^{*}|_{U_{2}})^{*}|_{U_{1}}=f\;\;\;\mbox{ for }f\in \Hom (U_{1},U_{2}^{*}).
\end{eqnarray}}
\er

For a $V\otimes V$-module $A$ and a $V$-module $W$, denote by ${\cal{M}}_{Q(z)}(A,W)$
the space of all $Q(z)$-intertwining maps from $A$ to $W$. Then we have:

\bt{tuniversal}
Let $A$ be a $V\otimes V$-module, $W$ a $V$-module and $f$ 
a $Q(z)$-intertwining map from 
$A$ to $\overline{W}$ $(= (W')^{*})$. Then $f^{*}$ restricted to 
$W'$ $(\subset ((W')^{*})^{*})$
is a $V$-homomorphism from $W'$ 
to $T_{Q(z)}'(A)$ $(\subset A^{*})$. Furthermore, 
the linear map $f\mapsto f^{*}|_{W'}$ 
is a linear isomorphism from ${\cal{M}}_{Q(z)}(A,W)$
onto $\Hom_{V}(W',T_{Q(z)}'(A))$. 
\et

\pf Let $f$ be a $Q(z)$-intertwining map from 
$A$ to $\overline{W}$. In view of Remark \ref{rclassicalfact},
$f^{*}|_{W'}$ is a linear map from $W'$ to $A^{*}$.

Let $v\in V,\; w'\in W',\;a\in A$. 
Let $k\in \N$ be such that $x^{k}Y(v,x_{0})w'\in W'[[x]]$.
Note that $k$ depends only on $v$ and $w'$.
By applying $\Res_{x_{0}}x_{0}^{k}$ to
(\ref{eqzjacobiA}) we get
\begin{eqnarray}
(x_{1}-z)^{k}\<f(Y_{1}^{o}(v,x_{1})a),w'\>=
(x_{1}-z)^{k}\<f(Y_{2}(v,x_{1})a),w'\>.
\end{eqnarray}
Then
\begin{eqnarray}
(x_{1}-z)^{k}\<f^{*}(w'),Y_{1}^{o}(v,x_{1})a\>=
(x_{1}-z)^{k}\<f^{*}(w'),Y_{2}(v,x_{1})a\>.
\end{eqnarray}
By Lemma \ref{lweakcommqzcomp}, $f^{*}(w')\in \tilde{T}_{Q(z)}'(A)$.
Therefore $f^{*}$ maps $W'$ into $\tilde{T}_{Q(z)}'(A)$.

Furthermore, using (\ref{eqzaction}) with $\lambda= f^{*}(w')$
and (\ref{eqzjacobiA}) we obtain
\begin{eqnarray}
& &\<Y_{Q(z)}'(v,x_{0})f^{*}(w'),a\>\nonumber\\
&=&\Res_{x_{1}}x_{0}^{-1}\delta\left(\frac{x_{1}-z}{x_{0}}\right)
\<f^{*}(w'),Y_{1}^{o}(v,x_{1})a\>
-\Res_{x_{1}}x_{0}^{-1}\delta\left(\frac{z-x_{1}}{-x_{0}}\right)
\<f^{*}(w'),Y_{2}(v,x_{1})a\>\nonumber\\
&=&\Res_{x_{1}}x_{0}^{-1}\delta\left(\frac{x_{1}-z}{x_{0}}\right)
\<f(Y_{1}^{o}(v,x_{1})a),w'\>
-\Res_{x_{1}}x_{0}^{-1}\delta\left(\frac{z-x_{1}}{-x_{0}}\right)
\<f(Y_{2}(v,x_{1})a),w'\>\nonumber\\
&=&\<Y^{o}(v,x_{0})f(a),w'\>\nonumber\\
&=&\<f(a),Y(v,x_{0})w'\>\nonumber\\
&=&\<f^{*}(Y(v,x_{0})w'),a\>.
\end{eqnarray}
Thus
\begin{eqnarray}
Y_{Q(z)}'(v,x_{0})f^{*}(w')=f^{*}(Y(v,x_{0})w')\;\;\;\mbox{ for }v\in
V,\; w'\in W'.
\end{eqnarray}
Therefore, $f^{*}|_{W'}$ is a $V$-homomorphism.
In view of Remark \ref{rclassicalfact}, the map $f\mapsto f^{*}|_{W'}$ is an injective
map from $\Hom_{Q(z)}(A,W)$ to $\Hom_{V}(W',T'_{Q(z)}(A))$.

On the other hand, let $g$ be a $V$-homomorphism from $W'$ to 
$\tilde{T}_{Q(z)}'(A)$ $(\subset A^{*})$.
Then we have a linear map $g^{*}|_{A}$ from
$A$ to $(W')^{*}=\overline{W}$.
Now we prove that $g^{*}|_{A}$ is a $Q(z)$-intertwining map.

Let $a\in A,\; v\in V,\; w'\in W'$. 
Using (\ref{e2.53}) with $\lambda=g(w')$ we get
\begin{eqnarray}
& &x_{0}^{-1}\delta\left(\frac{x_{1}-z}{x_{0}}\right)
\<g^{*}(Y_{1}^{o}(v,x_{1})a),w'\>
-x_{0}^{-1}\delta\left(\frac{z-x_{1}}{-x_{0}}\right)
\<g^{*}(Y_{2}(v,x_{1})a),w'\>\nonumber\\
&=&x_{0}^{-1}\delta\left(\frac{x_{1}-z}{x_{0}}\right)
\<g(w'),Y_{1}^{o}(v,x_{1})a\>
-x_{0}^{-1}\delta\left(\frac{z-x_{1}}{-x_{0}}\right)
\<g(w'),Y_{2}(v,x_{1})a\>\nonumber\\
&=&z^{-1}\delta\left(\frac{x_{1}-x_{0}}{z}\right)
\<Y_{Q(z)}'(v,x_{0})g(w'), a\>\nonumber\\
&=&z^{-1}\delta\left(\frac{x_{1}-x_{0}}{z}\right)
\<g(Y(v,x_{0})w'), a\>\nonumber\\
&=&z^{-1}\delta\left(\frac{x_{1}-x_{0}}{z}\right)
\<g^{*}(a),Y(v,x_{0})w'\>\nonumber\\
&=&z^{-1}\delta\left(\frac{x_{1}-x_{0}}{z}\right)
\<Y^{o}(v,x_{0})g^{*}(a),w'\>.
\end{eqnarray}
This proves that $g^{*}$ is a $Q(z)$-intertwining map from $A$ 
to $\overline{W}$. Then
the map $f\mapsto f^{*}|_{W'}$ is a linear isomorphism.
$\;\;\;\;\Box$

Combining Theorems \ref{tuniversal} and \ref{tqz} we immediately have:

\bc{cadjointproperty}
Let $A$ be a $V\otimes V$-module and $W$ a $V$-module.
Let $f$ be a $V\otimes V$-homomorphism from 
$A$ to $R_{Q(z)}(W')$ $(\subset (W')^{*})$. Then $f^{*}$ restricted to 
$W'$ $(\subset ((W')^{*})^{*})$
is a $V$-homomorphism from $W'$ 
to $T_{Q(z)}'(A)$ $(\subset A^{*})$. Furthermore, 
the linear map $f\mapsto f^{*}|_{W'}$ 
is a linear isomorphism from $\Hom_{V\otimes V}(A,R_{Q(z)}(W'))$
onto $\Hom_{V}(W',T_{Q(z)}'(A))$.$\;\;\;\;\Box$
\ec


We now assume that $V$ is rational in the sense of Huang and Lepowsky.
Then for each $V\otimes V$-module $A$, $T'_{Q(z)}(A)$ is a direct sum of
irreducible $V$-modules.
Denote by ${\cal{C}}^{o}_{V\otimes V}$ the category
of $V\otimes V$-modules on which both $L^{1}(0)$ and $L^{2}(0)$ 
act semisimply.
It follows immediately from ([Li2], Proposition 4.15) 
that any $V\otimes V$-module $A$ from
${\cal{C}}^{o}_{V\otimes V}$ is semisimple. 
Then each $V\otimes V$-module $A$ from ${\cal{C}}^{o}_{V\otimes V}$
is a direct sum of finitely many irreducible $V\otimes V$-modules.
Because by assumption all the fusion rules are finite, 
it follows from Corollary \ref{cintertwininghomo} that 
for any $V\otimes V$-module $A$ in ${\cal{C}}_{V\otimes V}^{o}$ 
and $V$-module $W$,
$\Hom _{Q(z)}(A,W)$ is finite-dimensional. Then
it follows from Theorem \ref{tuniversal} that 
the multiplicity of $W'$ in $T'_{Q(z)}(A)$ is finite.  
Therefore, $T'_{Q(z)}(A)$ is a $V$-module.
To summarize, we have (cf. [HL2], Proposition 5.6):

\bp{pmoduleT'}
Assume that $V$ is rational in the sense of Huang and Lepowsky.
Let $A$ be a $V\otimes V$-module in ${\cal{C}}^{o}_{V\otimes V}$.
Then $T'_{Q(z)}(A)$ is a $V$-module.$\;\;\;\;\Box$
\ep

Assume that $V$ is rational in the sense of Huang and Lepowsky.
Then we have a contravariant functor $T'_{Q(z)}$ from the category 
${\cal{C}}^{o}_{V\otimes V}$ of
$V\otimes V$-modules to the category of $V$-modules.
For a $V\otimes V$-module from ${\cal{C}}_{V\otimes V}^{o}$, we set
\begin{eqnarray}
T_{Q(z)}(A)=(T_{Q(z)}'(A))'.
\end{eqnarray}
This gives us a coinvariant functor $T_{Q(z)}$. 

\br{rhlhalftensor}
{\em In the case that $A=W_{1}\otimes W_{2}$ where $W_{1},W_{2}$ 
are $V$-modules, $T'_{Q(z)}(W_{1}\otimes W_{2})$ is exactly
the $V$-module $W_{1}\hboxtr_{Q(z)}W_{2}$ constructed in [HL2-3], 
so that $T_{Q(z)}(W_{1}\otimes W_{2})$ is exactly the Hunag-Lepowsky's
$Q(z)$-tensor product module of $W_{1}$ with $W_{2}$.}
\er

Note that $\Hom_{V}(W_{1},W_{2})$ is canonically isomorphic to
$\Hom_{V}(W_{2}',W_{1}')$ (see [HL2]).
Using Corollary \ref{cadjointproperty} we immediately have:

\bt{trightadjoint}
Assume that $V$ is rational in the sense of Huang-Lepowsky. Then
the functor ${\cal{F}}_{Q(z)}$ is a right adjoint of the functor
$T_{Q(z)}$. $\;\;\;\;\Box$
\et

Let $A$ be a $V\otimes V$-module and $W$ a $V$-module. 
{\em An intertwining operator} from $A$ to $W$ is
a linear map $F_{x}$ from $A$ to $W((x))$ such that the following
conditions hold for $v\in V,\; a\in A$:
\begin{eqnarray}
& &x_{0}^{-1}\delta\left(\frac{x_{1}-x_{2}}{x_{0}}\right)
Y(v,x_{1})F_{x_{2}}(a)
-x_{0}^{-1}\delta\left(\frac{x_{2}-x_{1}}{-x_{0}}\right)
F_{x_{2}}(Y_{2}(v,x_{1})a)\nonumber\\
&=&x_{2}^{-1}\delta\left(\frac{x_{1}-x_{0}}{x_{2}}\right)
F_{x_{2}}(Y_{1}(v,x_{0})a)
\end{eqnarray}
and
\begin{eqnarray}
F_{x}(L_{1}(-1)a)={d\over dx}F_{x}(a).
\end{eqnarray}
Denote by ${\cal{V}}(A,W)$ the space of all
intertwining operators from $A$ to $W$.

\br{rnewintermapsQ(z)}
{\em If $A=W_{1}\otimes W_{2}$ for $V$-modules $W_{1}$ and $W_{2}$,
then the same proof with necessary adjustments of [HL4] shows that 
${\cal{M}}_{Q(z)}(A,W)$ is canonically isomorphic to
${\cal{V}}(A,W)$ for any $V$-module $W$. 
But, for a general $V\otimes V$-module $A$ from the category 
${\cal{C}}_{V\otimes V}^{o}$, it is not easy to construct a
canonical linear isomorphism between ${\cal{V}}(A,W)$ and 
${\cal{M}}_{Q(z)}(A,W)$. }
\er


\section{Huang-Lepowsky's $P(z)$-tensor functor}
This section is parallel to Section 5. 
We here shall extend Huang-Lepowsky's $P(z)$-tensor functor 
constructed in [HL4] for the category of $V$-modules to 
a functor $T_{P(z)}$ from the category ${\cal{C}}_{V\otimes V}$ 
of $V\otimes V$-modules to the category of $V$-modules. 
In the construction of $T_{P(z)}$ we exploit a relation between
$T_{P(z)}$ and $T_{Q(-z^{-1})}$. 
The idea of this approach is due to [HL4], where the $V$-module
structure of the $P(z)$-tensor product module was established by using 
the $V$-module structure of the $Q(z^{-1})$-tensor product module.

Let $A$ be a weak $V\otimes V$-module for now, $A^{*}$ the dual space.
For $v\in V,\; \lambda\in A^{*}$ we define
$$Y_{P(z)}'(v,x)\lambda\in A^{*}[[x,x^{-1}]]$$
by
\begin{eqnarray}\label{epzaction}
& &\<Y_{P(z)}'(v,x)\lambda,a\>\nonumber\\
&=&\Res_{x_{0}}z^{-1}\delta\left(\frac{x^{-1}-x_{0}}{z}\right)
\<\lambda,Y_{1}(e^{xL(1)}(-x^{-2})^{L(0)}v,x_{0})a\>\nonumber\\
& &+ \Res_{x_{0}}x_{0}^{-1}\delta\left(\frac{z-x^{-1}}{-x_{0}}\right)
\<\lambda,Y_{2}^{o}(v,x)a\>\\
&=&\Res_{x_{0}}z^{-1}\delta\left(\frac{x^{-1}-x_{0}}{z}\right)
\<\lambda,Y_{1}(e^{xL(1)}(-x^{-2})^{L(0)}v,x_{0})a\>
+\<\lambda,Y_{2}^{o}(v,x)a\>\hspace{1cm}
\end{eqnarray}
(cf. [HL4], (13.6)) for $a\in A$.

\bd{dP(z)compatibility}
{\em A linear functional $\lambda$ on $A$ is said to satisfy
Huang-Lepowsky's {\em $P(z)$-compatibility condition} 
if the following conditions hold:

(a) The $P(z)$-{\em lower truncation condition}:
For all $v\in V$, the formal Laurent series $Y_{P(z)}'(v,x)\lambda$
involves only finitely many negative (integral) powers of $x$.

(b) For all $v\in V,\; a\in A$,
\begin{eqnarray}\label{ecompatibilityb}
& &x_{0}^{-1}\delta\left(\frac{x_{1}^{-1}-z}{x_{0}}\right)
\<Y_{P(z)}'(v,x_{1})\lambda,a\>\nonumber\\
&=&z^{-1}\delta\left(\frac{x_{1}^{-1}-x_{0}}{z}\right)
\<\lambda,Y_{1}(e^{x_{1}L(1)}(-x_{1}^{-2})^{L(0)}v,x_{0})a\>\nonumber\\
& &+x_{0}^{-1}\delta\left(\frac{z-x_{1}^{-1}}{-x_{0}}\right)
\<\lambda, Y_{2}^{o}(v,x_{1})a\>.
\end{eqnarray}}
\ed

Clearly, all the linear functionals on $A$ that satisfy
the Huang-Lepowsky's $P(z)$-compatibility condition 
form a subspace of $A^{*}$, which 
we denote by $\tilde{T}'_{P(z)}(A)$.

\br{ranotherform} {\em In (\ref{epzaction}),
by substituting $x$ and 
$v$ with $x^{-1}$ and $e^{xL(1)}(-x^{-2})^{L(0)}v$, respectively,
we obtain
\begin{eqnarray}\label{eanotherformaction}
& &\<(Y_{P(z)}')^{o}(v,x)\lambda,a\>\nonumber\\
&=&\Res_{x_{0}}z^{-1}\delta\left(\frac{x-x_{0}}{z}\right)
\<\lambda,Y_{1}(v,x_{0})a\>
+\Res_{x_{0}}x_{0}^{-1}\delta\left(\frac{z-x}{-x_{0}}\right)
\<\lambda, Y_{2}(v,x)a\>,
\end{eqnarray}
noting that $e^{-xL(1)}(-x^{2})^{L(0)}e^{xL(1)}(-x^{-2})^{L(0)}=1$ 
([FHL], (5.3.1)).
Similarly, in (\ref{ecompatibilityb}), by substituting $x_{1}$ and 
$v$ with $x_{1}^{-1}$ and $e^{x_{1}L(1)}(-x_{1}^{-2})^{L(0)}v$, 
respectively, we obtain
\begin{eqnarray}\label{eanotherform}
& &x_{0}^{-1}\delta\left(\frac{x_{1}-z}{x_{0}}\right)
\<(Y_{P(z)}')^{o}(v,x_{1})\lambda,a\>\nonumber\\
&=&z^{-1}\delta\left(\frac{x_{1}-x_{0}}{z}\right)
\<\lambda,Y_{1}(v,x_{0})a\>
+x_{0}^{-1}\delta\left(\frac{z-x_{1}}{-x_{0}}\right)
\<\lambda, Y_{2}(v,x_{1})a\>.
\end{eqnarray}
Clearly, this is also equivalent to (\ref{ecompatibilityb}).}
\er

In the following, we shall construct a weak $V\otimes V$-module 
$A_{new}$ from $A$ such that
\begin{eqnarray}
(\tilde{T}'_{P(z)}(A), Y_{P(z)}')
=(\tilde{T}'_{Q(-z^{-1})}(A_{new}), Y_{Q(-z^{-1})}').
\end{eqnarray}
Then it will follow immediately from Theorem \ref{texthl} that
the pair $(\tilde{T}'_{P(z)}(A), Y_{P(z)}')$
carries the structure of a weak $V$-module.
The idea of this approach is due to [HL4], where the $V$-module
structure of the $P(z)$-tensor product module was established by using 
the $V$-module structure of the $Q(z^{-1})$-tensor product module.

By replacing $x_{0}$ with $-zx_{0}x_{1}^{-1}$
in (\ref{ecompatibilityb}) we get
\begin{eqnarray}
& &-z^{-1}x_{1}x_{0}^{-1}\delta\left(\frac{-z^{-1}+x_{1}}{x_{0}}\right)
\<Y_{P(z)}'(v,x_{1})\lambda,a\>\nonumber\\
&=&z^{-1}\delta\left(\frac{z^{-1}+x_{0}}{x_{1}}\right)
\<\lambda,Y_{1}(e^{x_{1}L(1)}(-x_{1}^{-2})^{L(0)}v,-zx_{0}/x_{1})a\>
\nonumber\\
& &-z^{-1}x_{1}x_{0}^{-1}\delta\left(\frac{x_{1}-z^{-1}}{x_{0}}\right)
\<\lambda, Y_{2}^{o}(v,x_{1})a\>.
\end{eqnarray}
That is,
\begin{eqnarray}
& &x_{0}^{-1}\delta\left(\frac{-z^{-1}+x_{1}}{x_{0}}\right)
\<Y_{P(z)}'(v,x_{1})\lambda,a\>\nonumber\\
&=&x_{0}^{-1}\delta\left(\frac{x_{1}-z^{-1}}{x_{0}}\right)
\<\lambda, Y_{2}^{o}(v,x_{1})a\>\nonumber\\
& &-x_{1}^{-1}\delta\left(\frac{z^{-1}+x_{0}}{x_{1}}\right)
\<\lambda,Y_{1}(e^{x_{1}L(1)}(-x_{1}^{-2})^{L(0)}v,-zx_{0}/x_{1})a\>.
\end{eqnarray}
Using the fundamental properties of delta function we get
\begin{eqnarray}\label{ecompatibilityb1}
& &(-z^{-1})^{-1}\delta\left(\frac{x_{0}-x_{1}}{-z^{-1}}\right)
\<Y_{P(z)}'(v,x_{1})\lambda,a\>\nonumber\\
&=&x_{1}^{-1}\delta\left(\frac{x_{0}+z^{-1}}{x_{1}}\right)
\<\lambda, Y_{2}^{o}(v,x_{1})a\>\nonumber\\
& &-x_{1}^{-1}\delta\left(\frac{z^{-1}+x_{0}}{x_{1}}\right)
\<\lambda,Y_{1}(e^{x_{1}L(1)}(-x_{1}^{-2})^{L(0)}v,-zx_{0}/x_{1})a\>.
\end{eqnarray}
It is clear that (\ref{ecompatibilityb1}) is equivalent to
(\ref{ecompatibilityb}).
Furthermore, (\ref{ecompatibilityb1}) amounts to 
\begin{eqnarray}\label{ecompatibilityb3}
& &(-z^{-1})^{-1}\delta\left(\frac{x_{0}-x_{1}}{-z^{-1}}\right)
\<Y_{P(z)}'(v,x_{1})\lambda,a\>\nonumber\\
&=&x_{1}^{-1}\delta\left(\frac{x_{0}+z^{-1}}{x_{1}}\right)
\<\lambda, Y_{2}^{o}(v,x_{0}+z^{-1})a\>\nonumber\\
& &-x_{1}^{-1}\delta\left(\frac{z^{-1}+x_{0}}{x_{1}}\right)
\left\<\lambda,Y_{1}\left(e^{(z^{-1}+x_{0})L(1)}
(-(z^{-1}+x_{0})^{-2})^{L(0)}v,\frac{-zx_{0}}{z^{-1}+x_{0}}\right)a\right\>
\end{eqnarray}
(cf. (\ref{eqzcompatibility})).

The following results, formulated in [Li2] (Lemma 4.25, Remark 2.10), 
are special cases of the general result in change of
variable ([Z1-2], [H2]):

\bl{lgeneralL(0)}
Let $W$ be a weak $V$-module and $z$ a nonzero complex number. 
For $v\in V$, set
\begin{eqnarray}
Y[z](v,x)=Y(z^{L(0)}v,zx).
\end{eqnarray}
Then $W$ equipped with $Y[z]$ is a weak $V$-module.
Furthermore, if $W$ is a generalized $V$-module,
then $e^{L(0)\log z}$ 
is a $V$-isomorphism from $(W,Y_{W})$ to $(W,Y[z])$. 
\el

\bl{lgeneralL(1)}
Let $W$ be a weak $V$-module and $z$ a nonzero complex number. 
For $v\in V$, set
\begin{eqnarray}
Y\{z\}(v,x)=Y\left(e^{-z(1+zx)L(1)}(1+zx)^{-2L(0)}v,
\frac{x}{1+zx}\right).
\end{eqnarray}
Then $W$ equipped with $Y\{z\}$ is a weak $V$-module and
\begin{eqnarray}
(Y\{z\})^{o}(v,x)=Y^{o}(v,x+z)\;\left(=Y^{o}(v,y)|_{y=x+z}\right)
\;\;\;\mbox{ for }v\in V.
\end{eqnarray}
Furthermore, if $W$ is a weak $V$-module on which $L(1)$ acts semisimply,
then $e^{zL(1)}$ 
is a $V$-isomorphism from $(W,Y_{W})$ to $(W,Y\{z\})$. 
\el

Furthermore, we also have:

\bl{lgeneral2}
Let $W$ be a weak $V$-module and $z$ a nonzero complex number. 
For $v\in V$, set
\begin{eqnarray}
Y_{mix}(v,x)=Y\left(e^{(z^{-1}+x)L(1)}(z^{-1}+x)^{-2L(0)}(-1)^{L(0)}v,
\frac{-zx}{z^{-1}+x}\right).
\end{eqnarray}
Then $W$ equipped with $Y_{mix}$ is a weak $V$-module.
Furthermore, if $W$ is a generalized $V$-module on which $L(1)$ 
acts locally nilpotently, then $e^{(\pi i +\log z^{2})L(0)}e^{zL(1)}$ 
is a $V$-isomorphism from $(W,Y_{W})$ to $(W,Y_{mix})$. 
\el

\pf Recall from [HL4] (cf. [FHL], (5.3.1)-(5.3.3)) that on $V$,
\begin{eqnarray}
z_{0}^{L(0)}e^{xL(1)}=e^{z_{0}^{-1}xL(1)}z_{0}^{L(0)}
\end{eqnarray}
for any nonzero complex number $z_{0}$. Then,
for $v\in V$,
\begin{eqnarray}
Y_{mix}(v,x)&=&Y\left(e^{(z^{-1}+x)L(1)}(z^{-1}+x)^{-2L(0)}(-1)^{L(0)}v,
\frac{-zx}{z^{-1}+x}\right)\nonumber\\
&=&Y\left((-z^{2})^{L(0)}e^{-z(1+zx)L(1)}(1+zx)^{-2L(0)}v,
\frac{-z^{2}x}{1+zx}\right)\nonumber\\
&=&Y[-z^{2}]\left(e^{-z(1+zx)L(1)}(1+zx)^{-2L(0)}v,
\frac{x}{1+zx}\right)\nonumber\\
&=&(Y[-z^{2}])\{z\}(v,x),
\end{eqnarray}
where we are using the notions in Lemmas \ref{lgeneralL(0)} 
and \ref{lgeneralL(1)}.
Then it follows from Lemmas \ref{lgeneralL(0)} and \ref{lgeneralL(1)} 
that $(W,Y_{mix})$ is a weak $V$-module.

Assume that $W$ is a generalized $V$-module on which $L(1)$ 
acts locally nilpotently. 
Then
\begin{eqnarray}
Y_{mix}(v,x)=e^{L(0)\log (-z^{2})}
e^{zL(1)}Y(v,x)e^{-zL(1)}e^{-L(0)\log (-z^{2})}.
\end{eqnarray}
It is clear that the weak $V$-module structure $(W,Y_{mix})$ is 
the transported weak $V$-module structure through the linear map
$e^{L(0)\log (-z^{2})}e^{zL(1)}$.
$\;\;\;\;\Box$

Let $\theta$ be the natural automorphism of $V\otimes V$ uniquely
determined by
\begin{eqnarray}
\theta(u\otimes v)=v\otimes u\;\;\;\mbox{ for }u,v\in V.
\end{eqnarray}
Then we have a (weak) $V\otimes V$-module $A^{\theta}$, where
$A^{\theta}=A$ as a vector space and
\begin{eqnarray}
Y_{A^{\theta}}(a,x)=Y_{A}(\theta(a),x)\;\;\;\mbox{ for }a\in V\otimes V.
\end{eqnarray}

\bp{pnewbimodule}
Let $A$ be a weak $V\otimes V$-module and let $z$ be a nonzero complex number.
Let $Y_{new}$ be the linear map from $V\otimes V$ to $(\End A)[[x,x^{-1}]]$,
uniquely determined by
\begin{eqnarray}
Y_{new}(u\otimes v,x)
&=&Y_{2}\left(e^{(z^{-1}+x)L(1)}(z^{-1}+x)^{-2L(0)}(-1)^{L(0)}u,
\frac{-zx}{z^{-1}+x}\right) \nonumber\\
& &\cdot  
Y_{1}\left(e^{-z^{-1}(1+z^{-1}x)L(1)}(1+z^{-1}x)^{-2L(0)}(-1)^{L(0)}v,
\frac{x}{1+z^{-1}x}\right)
\end{eqnarray}
for $u,v\in V$. Then $(A,Y_{new})$ carries the structure of a weak 
$V\otimes V$-module.
Denote this weak $V\otimes V$-module by $A_{new}$.
Furthermore, if $A$ is a $V\otimes V$-module in the category
${\cal{C}}_{V\otimes V}^{o}$, then $A_{new}$ is a $V\otimes V$-module 
in ${\cal{C}}_{V\otimes V}^{o}$ and the linear map 
$e^{L(0)\log (-z^{2})}e^{zL(1)}\otimes e^{z^{-1}L(1)}$ 
is a $V\otimes V$-isomorphism from $A^{\theta}$ to $A_{new}$.
\ep

\pf For $v\in V$, set
\begin{eqnarray}
Y^{new}_{1}(v,x)
&=&Y_{2}\left(e^{(z^{-1}+x)L(1)}(z^{-1}+x)^{-2L(0)}(-1)^{L(0)}v,
\frac{-zx}{z^{-1}+x}\right),\\
Y^{new}_{2}(v,x)
&=&Y_{1}\left(e^{-z^{-1}(1+z^{-1}x)L(1)}(1+z^{-1}x)^{-2L(0)}(-1)^{L(0)}v,
\frac{x}{1+z^{-1}x}\right).
\end{eqnarray}
By Lemmas \ref{lgeneralL(1)} and \ref{lgeneral2}, 
$(A,Y_{1}^{new})$ and $(A,Y_{2}^{new})$ are weak $V$-modules.
Clearly, $Y_{1}^{new}$ still commutes with $Y_{2}^{new}$. Then,
the same proof of Proposition 4.6.1 of [FHL] shows
that $(A, Y_{new})$ is a weak $V\otimes V$-module.
The rest immediately follows from 
Lemmas \ref{lgeneralL(1)} and \ref{lgeneral2}.
$\;\;\;\;\Box$

In view of Lemma \ref{lgeneralL(1)} and Proposition \ref{pnewbimodule},
combining (\ref{ecompatibilityb3}) with (\ref{eqzcompatibility}) 
we immediately have:

\bp{pQ(z)=P(z)}
Let $A$ be a weak $V\otimes V$-module and $z$ a nonzero complex number. Then
\begin{eqnarray}
(\tilde{T}'_{P(z)}(A), Y'_{P(z)})=
(\tilde{T}'_{Q(-z^{-1})}(A_{new}), Y'_{Q(-z^{-1})}).\;\;\;\;\Box
\end{eqnarray}
\ep

Now, we have (cf. [HL4], Theorem 13.9):

\bt{tQ(z)=P(z)}
Let $A$ be a weak $V\otimes V$-module. Then 
$(\tilde{T}'_{P(z)}(A), Y'_{P(z)})$
carries the structure of a weak $V$-module. 
Furthermore, if $A$ is a $V\otimes V$-module 
in the category ${\cal{C}}_{V\otimes V}^{o}$, then
the linear map $\psi^{*}$ suitably restricted 
is a $V$-isomorphism from $(\tilde{T}'_{P(z)}(A), Y'_{P(z)})$
to $(\tilde{T}'_{Q(-z^{-1})}(A^{\theta}), Y'_{Q(-z^{-1})})$, where
$\psi$ is the linear automorphism of $A$ defined by
\begin{eqnarray}
\psi(a)=(e^{L(0)\log (-z^{2})}e^{zL(1)}\otimes e^{z^{-1}L(1)})a
\;\;\;\;\mbox{ for }a\in A
\end{eqnarray}
(cf. [HL4], (13.19)) and $\psi^{*}$ is the dual of $\psi$.
\et

\pf The first assertion immediately follows from 
Proposition \ref{pQ(z)=P(z)} and Theorem \ref{texthl}.
With Proposition \ref{pnewbimodule}, it follows from the functorial property that
the restriction of $\psi^{*}$ is a
$V$-isomorphism from $(\tilde{T}'_{Q(-z^{-1})}(A_{new}), Y'_{Q(-z^{-1})})$
to $(\tilde{T}'_{Q(-z^{-1})}(A^{\theta}), Y'_{Q(-z^{-1})})$.
Then it follows from Proposition \ref{pQ(z)=P(z)} that $\psi^{*}$ is a
$V$-isomorphism from $(\tilde{T}'_{P(z)}(A), Y'_{P(z)})$
to $(\tilde{T}'_{Q(-z^{-1})}(A^{\theta}), Y'_{Q(-z^{-1})})$.$\;\;\;\;\Box$

Let $A$ be a $V\otimes V$-module and $W$ a $V$-module. 
{\em A $P(z)$-intertwining map} from $A$ to $W$ is
a linear map $F$ from $A$ to $\overline{W}$ such that
for $v\in V,\; a\in A$:
\begin{eqnarray}
& &x_{0}^{-1}\delta\left(\frac{x_{1}-z}{x_{0}}\right)
Y(v,x_{1})F(a)
-x_{0}^{-1}\delta\left(\frac{z-x_{1}}{-x_{0}}\right)
F(Y_{2}(v,x_{1})a)\nonumber\\
&=&z^{-1}\delta\left(\frac{x_{1}-x_{0}}{z}\right)
F(Y_{1}(v,x_{0})a).
\end{eqnarray}
Clearly, all $P(z)$-intertwining maps from $A$ to $W$
form a subspace of $\Hom (A,\overline{W})$, which we 
denote by ${\cal{M}}_{P(z)}(A,W)$.

Let $F_{x}$ be an intertwining operator map from $A$ to $W$.
Clearly, the evaluation of $F_{x}$ at $x=e^{\log z}$ is a 
$P(z)$-intertwining map from $A$ to $W$. This gives rise to 
a linear map from ${\cal{V}}(A,W)$ to ${\cal{M}}_{P(z)}(A,W)$.

The same proof of Proposition 12.2 in [HL4] (cf. [HL2], Proposition 4.7)
gives:

\bp{pnewintermaps}
Let $A$ be a $V\otimes V$-module from the category 
${\cal{C}}_{V\otimes V}^{o}$ and $W$ a $V$-module. 
Then ${\cal{V}}(A,W)$ is canonically isomorphic to
${\cal{M}}_{P(z)}(A,W)$
through evaluation $x=e^{\log z}$ 
in intertwining operators. $\;\;\;\;\Box$
\ep

The same proof of Theorem 4.5 in [Li2] (cf. Theorem \ref{tqz}) gives:

\bp{pAuniversal}
Let $A$ be a $V\otimes V$-module and $W$ a $V$-module. 
Then a $P(z)$-intertwining map from $A$ to $\overline{W}$ 
exactly amounts to a $V\otimes V$-homomorphism from $A$ to 
$R_{P(z)}(W')\subset {\cal{D}}_{P(z)}(W')$. $\;\;\;\;\Box$
\ep

The following theorem is parallel to Theorem \ref{tuniversal}:

\bt{tintertwiningmap=Vhom}
Let $A$ be a weak $V\otimes V$-module, $W$ a $V$-module and
$F$ a linear map from $A$ to $\overline{W}$. Then $F$ is a
$P(z)$-intertwining map if and only if $F^{*}|_{W'}$ is a
$V$-homomorphism from $W'$ to $\tilde{T}'_{P(z)}(A)$ $(\subset A^{*})$.
Furthermore, the linear map $F\mapsto F^{*}|_{W'}$ 
is a linear isomorphism from ${\cal{M}}_{P(z)}(A,W)$
onto $\Hom_{V}(W',T_{P(z)}'(A))$. 
\et

\pf Assume that $F$ is a
$P(z)$-intertwining map. Then for $v\in V, \; a\in A, \; w'\in W'$,
\begin{eqnarray}
& &x_{0}^{-1}\delta\left(\frac{x_{1}-z}{x_{0}}\right)
\<Y(v,x_{1})F(a),w'\>
-x_{0}^{-1}\delta\left(\frac{z-x_{1}}{-x_{0}}\right)
\<F(Y_{2}(v,x_{1})a),w'\>\nonumber\\
&=&z^{-1}\delta\left(\frac{x_{1}-x_{0}}{z}\right)
\<F(Y_{1}(v,x_{0})a),w'\>.
\end{eqnarray}
That is,
\begin{eqnarray}\label{ethatis}
& &x_{0}^{-1}\delta\left(\frac{x_{1}-z}{x_{0}}\right)
\<F^{*}(Y^{o}(v,x_{1})w'),a\>
-x_{0}^{-1}\delta\left(\frac{z-x_{1}}{-x_{0}}\right)
\<F^{*}(w'),Y_{2}(v,x_{1})a\>\nonumber\\
&=&z^{-1}\delta\left(\frac{x_{1}-x_{0}}{z}\right)
\<F^{*}(w'),Y_{1}(v,x_{0})a)\>.
\end{eqnarray}
By taking $\Res_{x_{0}}$ we get
\begin{eqnarray}
& &\<F^{*}(Y^{o}(v,x_{1})w'),a\>\nonumber\\
&=&\<F^{*}(w'),Y_{2}(v,x_{1})a\>
+\Res_{x_{0}}z^{-1}\delta\left(\frac{x_{1}-x_{0}}{z}\right)
\<F^{*}(w'),Y_{1}(v,x_{0})a)\>.
\end{eqnarray}
Combining this with (\ref{eanotherformaction}) we obtain
\begin{eqnarray}\label{eequality}
\<(Y_{P(z)}')^{o}(v,x_{1})F^{*}(w'),a\>=\<F^{*}(Y^{o}(v,x_{1})w'),a\>,
\end{eqnarray}
which is equivalent to
\begin{eqnarray}\label{ehomomorphism}
Y_{P(z)}'(v,x_{1})F^{*}(w')=F^{*}(Y(v,x_{1})w').
\end{eqnarray}
Since $Y(v,x_{1})w'\in W'((x))$, it follows that $F^{*}(w')$
satisfies the $P(z)$-lower truncation condition. Furthermore,
(\ref{ethatis}) and (\ref{eequality}) imply (\ref{eanotherform}).
In view of Remark \ref{ranotherform}, $F^{*}(w')\in \tilde{T}'_{P(z)}(A)$.
{}From (\ref{ehomomorphism}), $F^{*}|_{W'}$ is a $V$-homomorphism.
In view of Remark \ref{rclassicalfact}, $F\mapsto F^{*}|_{W'}$ 
is an injective map from ${\cal{M}}_{P(z)}(A,W)$ to $\Hom_{V}(W',T_{P(z)}'(A))$.

Let $G$ be a $V$-homomorphism from $W'$ into
$\tilde{T}'_{P(z)}(A)$ $(\subset A^{*})$. For $v\in V,\; a\in A,\; w'\in W$, using 
(\ref{eanotherform}) we get
\begin{eqnarray}
& &x_{0}^{-1}\delta\left(\frac{x_{1}-z}{x_{0}}\right)
\<G(Y^{o}(v,x_{1})w'),a\>\nonumber\\
&=&x_{0}^{-1}\delta\left(\frac{z-x_{1}}{-x_{0}}\right)
\<G(w'),Y_{2}(v,x_{1})a\>
+z^{-1}\delta\left(\frac{x_{1}-x_{0}}{z}\right)
\<G(w'),Y_{1}(v,x_{0})a\>.
\end{eqnarray}
that is,
\begin{eqnarray}
& &x_{0}^{-1}\delta\left(\frac{x_{1}-z}{x_{0}}\right)
\<G^{*}(a),Y^{o}(v,x_{1})w'\>\nonumber\\
&=&x_{0}^{-1}\delta\left(\frac{z-x_{1}}{-x_{0}}\right)
\<G^{*}(Y_{2}(v,x_{1})a),w'\>
+z^{-1}\delta\left(\frac{x_{1}-x_{0}}{z}\right)
\<G^{*}(Y_{1}(v,x_{0})a),w'\>,
\end{eqnarray}
what is equivalent to,
\begin{eqnarray}
& &x_{0}^{-1}\delta\left(\frac{x_{1}-z}{x_{0}}\right)
\<Y(v,x_{1})G^{*}(a),w'\>\nonumber\\
&=&x_{0}^{-1}\delta\left(\frac{z-x_{1}}{-x_{0}}\right)
\<G^{*}(Y_{2}(v,x_{1})a),w'\>
+z^{-1}\delta\left(\frac{x_{1}-x_{0}}{z}\right)
\<G^{*}(Y_{1}(v,x_{0})a),w'\>.
\end{eqnarray}
That is, $G^{*}|_{A}$ is a $P(z)$-intertwining map from $A$ to $\overline{W}$.
By Remark \ref{rclassicalfact}, we have $(G^{*}|_{A})^{*}|_{W'}=G$. Therefore,
$F\mapsto F^{*}|_{W'}$ 
is a linear isomorphism from ${\cal{M}}_{P(z)}(A,W)$ onto $\Hom_{V}(W',T_{P(z)}'(A))$.
$\;\;\;\;\Box$

Combining Proposition \ref{pAuniversal} with 
Theorem \ref{tintertwiningmap=Vhom} we immediately have 
(cf. Corollary \ref{cadjointproperty}):

\bc{cpreP(z)adjoint}
Let $A$ be a $V\otimes V$-module and $W$ a $V$-module. 
Let $f$ be a $V\otimes V$-homomorphism from $A$ to $R_{P(z)}(W')$ $(\subset (W')^{*})$. 
Then $f^{*}$ restricted to $W'$ $(\subset ((W')^{*})^{*})$ is a $V$-homomorphism from $W'$ 
to $T'_{P(z)}(A)$ $(\subset A^{*})$. Furthermore, 
the linear map $f\mapsto f^{*}|_{W'}$ is a linear isomorphism from
$\Hom_{V\otimes V}(A,R_{P(z)}(W'))$on to $\Hom_{V}(W',T'_{P(z)}(A))$.
$\;\;\;\;\Box$
\ec

With Theorem \ref{tintertwiningmap=Vhom}, using the same arguments 
of Proposition \ref{pmoduleT'} we obtain:

\bp{pt'p(z)module}
Assume that $V$ is rational in the sense  of Huang-Lepowsky.
Then for any $V\otimes V$-module $A$ in ${\cal{C}}_{V\otimes V}^{o}$, 
$T'_{P(z)}(A)$ is an (ordinary) $V$-module.
$\;\;\;\;\Box$
\ep

Assume that $V$ is rational in the sense  of Huang-Lepowsky.
For a $V\otimes V$-module $A$ in ${\cal{C}}_{V\otimes V}^{o}$, we set
\begin{eqnarray}
T_{P(z)}(A)=(T'_{P(z)}(A))'.
\end{eqnarray}
This gives us a coinvariant functor $T_{P(z)}$ from 
${\cal{C}}_{V\otimes V}^{o}$ to ${\cal{C}}_{V}$.

Immediately from Corollary \ref{cpreP(z)adjoint} we have:

\bt{trightadjointP(z)}
Assume that $V$ is rational in the sense  of Huang-Lepowsky.
Then the functor $T_{P(z)}$ is a right adjoint of the functor
${\cal{F}}_{P(z)}$. $\;\;\;\;\Box$
\et

\end{document}